\RequirePackage{fix-cm}

\documentclass[11pt,letterpaper]{article}
\newcommand{\argmin}{\operatornamewithlimits{arg\ min}}
\usepackage{graphicx,color}
\usepackage{bm}
\usepackage[lofdepth,lotdepth]{subfig}
\usepackage{hyperref}
\usepackage{amsmath, amsthm, amssymb}
\makeatletter
\def\@endtheorem{\endtrivlist}
\makeatother
\newtheorem{theorem}{Theorem}[section]
\usepackage{hyperref}
\newtheorem{lemma}[theorem]{Lemma}

\newtheorem{corollary}[theorem]{Corollary}
\newtheorem{definition}[theorem]{Definition}
\def\@endtheorem{\endtrivlist}
\begin{document}

\title{Convergence Rate for the Ordered Upwind Method\thanks{This work has been supported by the Ontario Government, and Natural Sciences and Engineering Research Council. The final publication is available at Springer via \url{http://dx.doi.org/10.1007/s10915-016-0163-3}}
}


%
\author{Alex Shum\thanks{Department of Applied Mathematics, \href{mailto:a5shum@uwaterloo.ca}{a5shum@uwaterloo.ca}, \href{mailto:kmorris@uwaterloo.ca}{kmorris@uwaterloo.ca}} \and
        Kirsten Morris\footnotemark[2] \and
				Amir Khajepour \thanks{Department of Mechanical and Mechatronics Engineering, \href{mailto:akhajepour@uwaterloo.ca}{akhajepour@uwaterloo.ca}, University of Waterloo, Waterloo, Canada, N2L 3G1} 
}

\date{}

%
%
\maketitle

\begin{abstract}
The Ordered Upwind Method (OUM) is used to approximate the viscosity solution of the static Hamilton-Jacobi-Bellman (HJB) with direction-dependent weights on unstructured meshes. The method has been previously shown to provide a solution that converges to the exact solution, but no convergence rate has been theoretically proven. In this paper, it is shown that the solutions produced by the OUM in the boundary value formulation converge at a rate of at least the square root of the largest edge length in the mesh in terms of maximum error. An example with similar order of numerical convergence is provided.
\end{abstract}
\section{Introduction}

The static Hamilton-Jacobi-Bellman (HJB) equation with a prescribed value on the boundary of a region $\Omega \subset \mathbb{R}^n$ where the solution is found on the interior of $\Omega$ arises in a number of optimization problems. Applications include optimal escape from a region \cite{kaltonthesis}, area patrol and perimeter surveillance \cite{oumboundarytracking}, modelling folds in structural geology \cite{oumfoldsinstructuralgeology} and reactive fluxes \cite{oumreactivefluxes}. 

There are two classes of semi-Lagrangian approximations \cite{kumarvlad} that approximate a solution to the static HJB equation. These approximations are known as semi-Lagrangian because the solution is approximated along short segments of characteristics dependent on the discretization. Both are solved on a fixed simplicial mesh or grid that discretizes the region of interest. The difference between them is the method in which the control is approximated.

In the first approach, the control is assumed to be held constant within an element of a mesh \cite{gonrof}. Non-iterative schemes such as the Ordered Upwind Method (OUM), Monotone Acceptance Ordered Upwind Method (MAOUM) \cite{kaltonpaper} and Fast Marching Method (FMM) \cite{fmmtriangle} use this approximation. In OUM, MAOUM and FMM, the order in which the solution on the vertices of the mesh (or grid) is found explicitly much like in Dijkstra's algorithm \cite{dijkstra} resulting in a significant speed up in computation, despite the coupling between vertices.

In the other semi-Lagrangian approximation, the control is assumed to be held fixed for a small time $\triangle t$. To determine the solution at a mesh point, a first-order reconstruction from nearby points on the discretization is required. An error bound $\mathcal{O}(\triangle t)$ has been shown for controls that have bounded variation \cite{bardifalconeerror}. Results of higher-order convergence rates using higher-order semi-Lagrangian approximation schemes of this type exist \cite{falconeferretti}. Many iterative algorithms \cite{bardi,bfm} have been devised that use this approximation.

Convergence rate results exist for the related time-dependent Hamilton-Jacobi equation, where similar half-order convergence is observed in terms of the longest time step (rather than edge length). These results have been proven for grid like discretizations \cite{crandalllionstwoapprox,souganidis} and have been extended to the use of triangular meshes \cite{abgrall} both using finite difference schemes. In \cite{bardi}, convergence rate results are given using similar schemes that include both time step and spatial discretizations. The proof of the main result in this work draws on some similar ideas such as doubling the variables in the use of an auxiliary function as in \cite{bardi} and \cite[Chapter 10]{evans}.

It is proven in this paper that the convergence rate of the approximate solution provided by OUM to the viscosity solution of the static HJB boundary value problem is at least $\mathcal{O}(\sqrt{h_{max}})$  in terms of maximum error, where $h_{max}$ is the longest edge length of a mesh. In \cite{oum}, the OUM was shown to provide an approximate solution to the static HJB equation that converges as $h_{max} \rightarrow 0$, but no convergence rate was obtained. The proof in this work is based on a similar result for FMM in \cite{monneau}. The OUM however is a different algorithm used to solve a wider class of problems where the weight (or speed) function can depend on position and direction and the boundary function can depend on position. The result in \cite{monneau} is proven on a uniform grid whereas the result here holds on a simplicial mesh. Simplicial meshes are better suited towards discretizing regions with complex geometries. A finer discretization may be required to obtain the same accuracy when the discretization is restricted to grids. A key step in the proof for the OUM convergence rate is showing the existence of a directionally complete stencil that is consistent with the result of OUM, an idea which was first presented in \cite{kaltonpaper}.

The optimal control problem along with an introduction to viscosity solutions will be presented in section \ref{secoptcontvisco}. In section \ref{secmeshes}, a general discretization of $\Omega \subset \mathbb{R}^n$, known as a simplicial mesh, will be described. The Ordered Upwind Method \cite{oum} will be reviewed in section \ref{secoum}. Properties of the OUM algorithm required in the proof of the main result will be presented in section \ref{secoumprop}. The convergence rate result will be proven in section \ref{secerrorbound}. An example demonstrating numerical convergence close to the proven theoretical rate will be presented in section 7. Conclusions and directions of future work will be discussed in section 8.

\section{Problem Formulation} \label{secoptcontvisco}
A point is denoted $\textbf{x} \in \mathbb{R}^n$ and the Euclidean norm is denoted $\left\|\cdot\right\|$. The set of positive real numbers is denoted $\mathbb{R}_+$. 
Let $\Omega \subset \mathbb{R}^n$ be open, connected, bounded with non-empty interior and boundary $\partial \Omega$. Let $\overline{\Omega} = \Omega \cup \partial \Omega$ be the closure of $\Omega$. 

Let $\mathcal{U} = \{\textbf{u}(\cdot): \mathbb{R}_+\cup\{0\} \rightarrow \mathbb{S}^{n-1}|\textbf{u}(\cdot) \text{ is measurable}\}$ where $\mathbb{S}^{n-1} = \{\textbf{u} \in \mathbb{R}^n |$ $\left\|\textbf{u}\right\|= 1 \}$ be the set of admissible controls and the trajectory $\textbf{y}: \mathbb{R}_+\cup\{0\} \rightarrow \overline{\Omega}$ is governed by control $\textbf{u}(\cdot) \in \mathcal{U}$,
\begin{equation} \dot{\textbf{y}}(t) = \textbf{u}(t), \textbf{y}(0) = \textbf{x}_0, \ \ \textbf{x}_0 \in \overline{\Omega}. \label{vehicleode} \end{equation}
The control problem is to steer $\textbf{y}(\cdot)$ from $\textbf{x}_0 \in \overline{\Omega}$ to any point on the boundary $\textbf{x}_f \in \partial \Omega$. The trajectory with initial condition $\textbf{y}(0) = \textbf{x}_0$ may be written $\textbf{y}_{\textbf{x}_0}(\cdot)$.

\begin{definition}
The \textbf{exit-time} $T: \overline{\Omega} \times \mathcal{U} \rightarrow \mathbb{R}_+\cup\{0\}$ is the first time $\textbf{y}_{\textbf{x}_0}(\cdot)$ reaches $\textbf{x}_f \in \partial\Omega$ under the influence of the control $\textbf{u}(\cdot)$,
\begin{equation} T(\textbf{x}_0, \textbf{u}(\cdot)) = \inf \{t | \textbf{y}_{\textbf{x}_0}(t) \in \partial \Omega \}. \end{equation} \end{definition}

To discuss optimality, a cost is assigned to each control.
\begin{definition}
The \textbf{cost function}, Cost: $ \overline{\Omega} \times \mathcal{U} \rightarrow \mathbb{R}$ is
\begin{equation}
\text{Cost}(\textbf{x}_0,\textbf{u}(\cdot)) = \int_0^{T(\textbf{x}_0,\textbf{u}(\cdot))} g(\textbf{y}_{\textbf{x}_0}(s),\textbf{u}(s))ds  + q(\textbf{y}_{\textbf{x}_0}(T(\textbf{x}_0,\textbf{u}(\cdot)))), \text{ for } \textbf{x}_0 \in \overline{\Omega}
\label{vehiclecost}
\end{equation}
 where $q: \partial\Omega \rightarrow \mathbb{R}$ is the boundary exit-cost and $g: \overline{\Omega} \times \mathbb{S}^{n-1} \rightarrow \mathbb{R}_+$ is the weight.
\end{definition}
 The optimal control problem is to find a control $\textbf{u}^*(\cdot)$ that minimizes (\ref{vehiclecost}).
\begin{definition}
The \textbf{value function} $V: \overline{\Omega} \rightarrow \mathbb{R}$ at $\textbf{x} \in \overline{\Omega}$ is the cost associated with the optimal control  $\textbf{u}^*(\cdot)$ for reaching any $\textbf{x}_f \in \partial \Omega$ from $\textbf{x}$,
\begin{equation} V(\textbf{x}) = \inf_{\textbf{u}(\cdot) \in \mathcal{U}} \text{Cost}(\textbf{x},\textbf{u}(\cdot)).
\label{vehiclevalue} \end{equation}
\end{definition}
The value function at $\textbf{x} \in \overline{\Omega}$ is the lowest cost to reach $\partial \Omega$ from $\textbf{x}$. The value function satisfies the  continuous {\em Dynamic Programming Principle} (DPP).

\begin{theorem} \label{bellmancont} (Dynamic Programming Principle \cite[Theorem 10.3.1]{evans})
For $h > 0$, $t \geq 0$, such that $0 \leq t+h \leq T(\textbf{x}_0,\textbf{u}^*(\cdot))$, 
\begin{equation} \displaystyle 
V(\textbf{y}_{\textbf{x}_0}(t)) = \inf_{\textbf{u}(\cdot) \in \mathcal{U}} \left\{ \int_t^{t+h} g(\textbf{y}_{\textbf{x}_0}(s),\textbf{u}(s)) ds + V(\textbf{y}_{\textbf{x}_0}(t+h)) \right\}. 
\label{contopt}
\end{equation}
\end{theorem}

For $V$ to be continuous on $\overline{\Omega}$, continuity between $V$ on $\Omega$ and $q$ on $\partial \Omega$ must be established. Let $L: \overline{\Omega} \times \overline{\Omega}$ be
\begin{equation} \label{Lequation} \displaystyle L(\textbf{x}_1,\textbf{x}_2) = \inf_{\textbf{u}(\cdot) \in \mathcal{U}} \left\{ \int_0^{\tau} g(\textbf{y}_{\textbf{x}_1}(s),\textbf{u}(s))ds \text{ } \Big | \text{ }\textbf{y}_{\textbf{x}_1}(\tau) = \textbf{x}_2,  \textbf{y}_{\textbf{x}_1}(t) \in \overline{\Omega}, t \in (0,\tau) \right\}. \end{equation}

\begin{definition} \label{compatconddef}
The exit-cost $q$ is \textbf{compatible} (with the continuity of $V$) if

\begin{equation} q(\textbf{x}_1) - q(\textbf{x}_2) \leq L(\textbf{x}_1, \textbf{x}_2) \label{compatcond} \end{equation}
for all $\textbf{x}_1,\textbf{x}_2 \in \partial \Omega$.
\end{definition}

\begin{definition} \label{speedprofiledef}
The \textbf{speed profile} of $g(\textbf{x},\textbf{u})$ is
$$\displaystyle \mathcal{U}_g(\textbf{x}) = \left\{\frac{t\textbf{u}}{g(\textbf{x},\textbf{u})} \Big | \textbf{u} \in \mathbb{S}^{n-1} \text{ and } t \in [0,1] \right\}.$$
\end{definition}
In $\mathbb{R}^2$, the speed profile is the shape centred at $\textbf{x}$ with radius $1/g(\textbf{x},\textbf{u})$ at the angle corresponding to the direction $\textbf{u}$. 

The optimal control problem (\ref{vehicleode}), (\ref{vehiclecost}) will be assumed to satisfy the following:
\begin{description}
\item{\textbf{(P1)} The boundary function $q$ is compatible with the continuity of $V$.}

\item{\textbf{(P2)} There exist constants $G_{min}, G_{max} \in \mathbb{R}_+$ and continuous functions $g_{min}, g_{max}: \overline{\Omega} \rightarrow \mathbb{R}_+$ such that for all $\textbf{x} \in \overline{\Omega}$ and $\textbf{u} \in \mathbb{S}^{n-1}$,
\begin{equation} \label{gbounds} 0 < G_{min} \leq g_{min}(\textbf{x}) \leq g(\textbf{x}, \textbf{u}) \leq g_{max}(\textbf{x}) \leq G_{max} < \infty. \end{equation}}
\item{\textbf{(P3)} There exists $L_g \in \mathbb{R}_+$ such that for $\textbf{x}_1, \textbf{x}_2 \in \overline{\Omega}$ and $\textbf{u} \in \mathbb{S}^{n-1}$,
\begin{equation} \label{lipschitzweight}
|g(\textbf{x}_1,\textbf{u}) - g(\textbf{x}_2,\textbf{u})| \leq L_g\left\|\textbf{x}_1 - \textbf{x}_2\right\|.
\end{equation}
}
\item{\textbf{(P4)} For all $\textbf{x}_1,\textbf{x}_2 \in \overline{\Omega}$ and $\lambda \in (0,1)$, $\lambda\textbf{x}_1 + (1-\lambda)\textbf{x}_2  \in \overline{\Omega}$. }

\item{\textbf{(P5)}  The speed profile $\mathcal{U}_g(\textbf{x})$ is convex for all $\textbf{x} \in \Omega$. 

Assumption \textbf{(P5)} is needed to guarantee uniqueness in the optimizing direction in the approximated problem provided $\nabla V$ exists \cite{kaltonpaper,vladthesis}.}
\end{description}

\begin{lemma}
The boundary function $q: \partial \Omega \rightarrow \mathbb{R}$ is Lipschitz-continuous. \label{qlipschitz}
\end{lemma}
The proof follows from \textbf{(P1)},\textbf{(P2)}, and \textbf{(P4)} with Lipschitz constant $2G_{max}$.

Since $q$ is Lipschitz-continuous on a compact subset of $\mathbb{R}^n$, there exist $q_{min}, q_{max} \in \mathbb{R}$  such that
\begin{equation}\label{qbound} q_{min} \leq q(\textbf{x}) \leq q_{max}. \end{equation}

Define the Hamiltonian $H: \Omega \times \mathbb{R}^n \rightarrow \mathbb{R}$
\begin{equation} \displaystyle
H(\textbf{x}, \textbf{p}) = -\min_{\textbf{u} \in \mathbb{S}^{n-1}} \{ \textbf{p} \cdot \textbf{u} + g(\textbf{x},\textbf{u}) \}.
\label{hamdef}
\end{equation}

\noindent The corresponding static Hamilton-Jacobi-Bellman (HJB) equation which can be derived from a first-order approximation of (\ref{contopt}) \cite{vladthesis} is
\begin{equation}  H(\textbf{x},\nabla V)=\min_{\textbf{u} \in \mathbb{S}^{n-1}} \{(\nabla V(\textbf{x}) \cdot \textbf{u})+ g(\textbf{x},\textbf{u})\} = 0, \textbf{x} \in \Omega,  \label{conttheo1} \end{equation}
\begin{equation*} V(\textbf{x}) = q(\textbf{x}), \text{ for } \textbf{x} \in \partial \Omega. \end{equation*}

\begin{definition} \label{chardirdef}
The \textbf{characteristic direction} $\textbf{u}^*: \Omega \rightarrow \mathbb{S}^{n-1}$ at $\textbf{x} \in \Omega$ is an optimizer of (\ref{conttheo1}) at \textbf{x}.
\end{definition}

Even for smooth $g(\textbf{x},\textbf{u})$, $q(\textbf{x})$ and $\partial \Omega$, $\nabla V$ (and hence unique $\textbf{u}^*$) may not exist over all of $\Omega$.  The weaker notion of viscosity solutions \cite{bardi}, is used to describe solutions of (\ref{hamdef}). Let $C^k(\Omega)$, $k\in\mathbb{N}\cup\{\infty\}$ denote the space of functions on $\Omega$ that are $k$-times continuously-differentiable.
\begin{definition} \cite{bardi} \label{subsol}
A function $\underline{V}: \overline{\Omega} \rightarrow \mathbb{R}$ is a \textbf{viscosity subsolution} of (\ref{conttheo1}) if for any $\phi \in C^\infty(\Omega)$, 
\begin{equation} \label{subsoleq} H(\textbf{x}_0, \nabla \phi(\textbf{x}_0)) \leq 0, \end{equation}
at any local maximum point $\textbf{x}_0 \in \Omega$ of $\underline{V} - \phi$.
\end{definition}
\begin{definition} \cite{bardi} \label{supersol}
A function $\overline{V}: \overline{\Omega} \rightarrow \mathbb{R}$ is a \textbf{viscosity supersolution} of (\ref{conttheo1}) if for any $\phi \in C^\infty(\Omega)$, 
\begin{equation} \label{supersoleq} H(\textbf{x}_0, \nabla \phi(\textbf{x}_0)) \geq 0, \end{equation}
at any local minimum point $\textbf{x}_0 \in \Omega$ of $\overline{V} - \phi$.
\end{definition}
\begin{definition} \cite{bardi}
A \textbf{viscosity solution} of the static HJB (\ref{conttheo1}) is both a viscosity subsolution and a viscosity supersolution of (\ref{conttheo1}).
\label{viscodef}
\end{definition}

\section{Simplicial Meshes} \label{secmeshes}
Viscosity solutions are often difficult to find analytically. The region $\overline{\Omega}$ will be discretized using a simplicial mesh on which $V$ (\ref{vehiclevalue}) will be solved approximately. 
 
\begin{definition}
A set of points $F = \{\textbf{x}_0,...,\textbf{x}_k \} \subset \mathbb{R}^n$ is \textbf{affinely independent} if the vectors $\{\textbf{x}_1 - \textbf{x}_0$, ... , $\textbf{x}_k - \textbf{x}_0 \}$ are linearly independent. 
\end{definition}
\begin{definition}
A \textbf{$k$-simplex} (plural $k$-simplices) $\textbf{s} = \textbf{x}_0^{\textbf{s}}\textbf{x}_1^{\textbf{s}}\cdots\textbf{x}_k^{\textbf{s}}$ is the convex hull of an affinely independent set of points $F=\{\textbf{x}_0^{\textbf{s}},\textbf{x}_1^{\textbf{s}}...,\textbf{x}_k^{\textbf{s}} \}$.
\end{definition}


\begin{definition}
Suppose \textbf{s} is a $k$-simplex defined by the convex hull of $F$. A \textbf{face} of \textbf{s} is any $m$-simplex ($-1\leq m \leq k$) forming the convex hull of a subset of $F$ containing $m+1$ elements.
\end{definition}

\begin{definition}
A \textbf{simplicial mesh}, $X$ is a set of simplices such that
\begin{enumerate}
\item Any face of a simplex in $X$ is also in $X$.
\item The intersection of two simplices $\textbf{s}_1,\textbf{s}_2 \in X$ is a face of $X$.
\end{enumerate}
\end{definition}
\begin{definition}
A \textbf{$k$-simplicial mesh} is a simplicial mesh where the highest dimension of any simplex in $X$ is $k$.
\end{definition}
Denote $X_j$, $0 \leq j \leq n$ the set of $j$-simplices of $X$. Elements of $X_0$, the $0$-simplices of $X$ are denoted $\textbf{x}_i$ and known as vertices. Elements of $X_1$, the $1$-simplices of $X$, are known as edges. 

Suppose $X \subset \mathbb{R}^n$ is an $n$-simplicial mesh. For $0 \leq k \leq n$, define
\begin{equation} \label{baryset} \displaystyle \Xi_k = \left\{ (\zeta_0,\zeta_1,...,\zeta_k) \in \mathbb{R}^{k+1} \Big| \sum_{j=0}^k \zeta_j = 1, \zeta_j \in [0,1] \text{ } \forall \text{ } 0 \leq j \leq {k-1} \right\}. \end{equation}
\begin{definition}
The \textbf{barycentric coordinates} of $\textbf{x}\in \mathbb{R}^n$ belonging to a $k$-simplex $\textbf{s}$ is a vector $\zeta = (\zeta_0,...,\zeta_{k}) \in \Xi_k$ such that $\textbf{x} = \sum_{j=0}^{k} \zeta_j\textbf{x}_j^{\textbf{s}}$.
\end{definition}

\begin{definition} \label{contained}
A closed region $A \subset \mathbb{R}^n$ is \textbf{contained} in an $n$-simplicial mesh $X$ if for every $\textbf{x} \in A$, there exists $\textbf{s} = \textbf{x}_0^{\textbf{s}}\textbf{x}_1^{\textbf{s}}\cdots\textbf{x}_n^{\textbf{s}}$ and $\zeta = (\zeta_0,\zeta_1,...,\zeta_n) \in \Xi_{n}$ such that $\textbf{x} = \sum_{j=0}^n \xi_j \textbf{x}_j^{\textbf{s}}$. 
\end{definition}
\begin{definition} \label{hmaxdef}
The \textbf{maximum edge length} $h_{max}$ is the length of the longest edge of $X$.
\end{definition}

\begin{definition} \label{hmindef}
{Let $1 \leq k\leq n$. A \textbf{neighbour} of simplex $\textbf{x}_0\textbf{x}_1\cdots\textbf{x}_{k-1} \in X_{k-1}$, is a vertex $\textbf{x}_k \in X_0$ such that $\textbf{x}_0\textbf{x}_1\cdots\textbf{x}_k \in X_k$.}
\end{definition}

\begin{definition}
{The \textbf{minimum simplex height} $h_{min}$ of $X$ is the shortest perpendicular distance between any $\textbf{s} \in X_{n-1}$ with its neighbours.}
\end{definition}

\begin{figure}
\begin{center} \includegraphics[height = 1.5in]{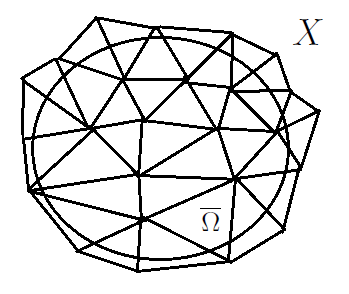} \end{center}
\caption{An example of $\overline{\Omega} \subset \mathbb{R}^2$ contained in a 2-simplicial mesh $X$.}
\label{meshcover}
\end{figure}

If $n=2$, then $h_{min}$ is the shortest triangle height. The following assumptions will be made on the ($n$-simplicial) mesh $X \subset \mathbb{R}^n$ on which the approximation of $V$ in the optimal control problem (\ref{vehicleode}), (\ref{vehiclecost}) will be found.   
\begin{description}
\item \textbf{(M1)} There exists $M \in \mathbb{R}_+$ such that $1 \leq \frac{h_{max}}{h_{min}} \leq M$.
\item \textbf{(M2)} The region $\overline{\Omega}$ is contained (Definition \ref{contained}) in the mesh $X$.
\item \textbf{(M3)} The mesh $X$ is bounded and has a finite number of vertices $X_0$.
\end{description}
The value $M$ is a measure of the worst-case degeneracy for a mesh $X$. An example of $\overline{\Omega} \subset \mathbb{R}^2$ being contained in a mesh $X$ is shown in Figure \ref{meshcover}. With the discretization definitions and assumptions stated, the OUM will now be presented.
\section{Review of the Ordered Upwind Method} \label{secoum}
The OUM \cite{oum} is used to find an approximation $\widetilde{V}: X_0 \rightarrow \mathbb{R}$ of $V$ in (\ref{contopt}) on the vertices of an $n$-simplicial mesh $X \subset \mathbb{R}^n$ satisfying \textbf{(M1)} -\textbf{(M3)}. 

The vertices of $X_0$ are assigned and updated between the following labels throughout the execution of the OUM. 
\begin{description}
\item \textit{\textbf{Far}} -  These vertices have values $\widetilde{V}(\textbf{x}_i) =  K$, where $K$ is a large value.  Computation of $\widetilde{V}$ has not yet started.
\item \textit{\textbf{Considered}} - These vertices have tentative values $\widetilde{V} < K$ and are computed using an update formula. 
\item \textit{\textbf{Accepted}} - These vertices have finalized values $\widetilde{V}$.
\end{description}
At any instant of the algorithm, each vertex in $X$ must be labelled exactly one of \textit{Accepted}, \textit{Considered} or \textit{Far}. Simplices with \textit{Accepted} label are further classified. 

\begin{figure}
\begin{center} \includegraphics[width=2in]{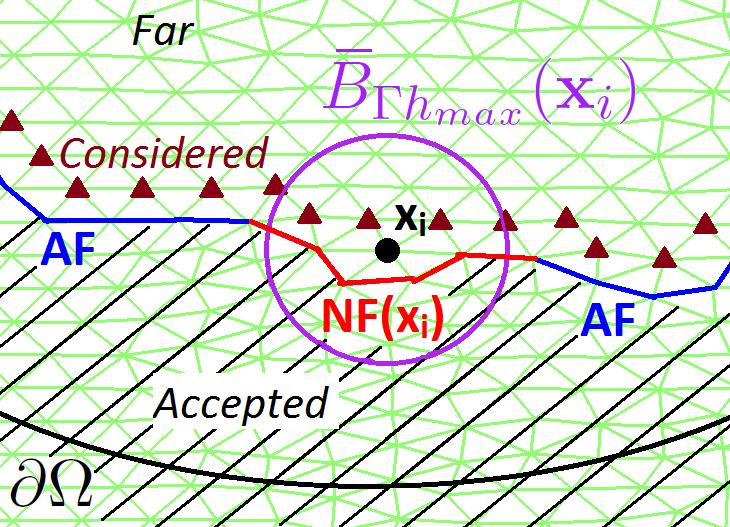} \end{center}
\caption{OUM Labels - An example for $\overline{\Omega} \subset \mathbb{R}^2$. The vertex $\textbf{x}_i$ with \textit{Considered} label is updated from the set of directions provided by $\textbf{NF}(\textbf{x}_i)$. Vertices labelled \textit{Accepted} are shaded, including vertices on the edges that make up \textbf{AF} and the \textit{Near Front} of $\textbf{x}_i$, $\textbf{NF}(\textbf{x}_i)$.  Vertices outside $\Omega$ are also labelled \textit{Accepted}. $\overline{B}_{\Gamma h_{max}}(\textbf{x}_i)$ is the closed ball with radius $\Gamma h_{max}$ and centre $\textbf{x}_i$. Vertices labelled \textit{Considered} are marked with a triangle. Unmarked vertices are labelled \textit{Far}.}
\label{classify}
\end{figure}

\begin{description}
\item \textit{\textbf{Accepted Front}} - The subset of vertices $X_0$ with \textit{Accepted} label that have a neighbour labelled \textit{Considered}. 
\item \textit{\textbf{AF}} - The subset of $X_{n-1}$ made of vertices on the \textit{Accepted Front} that have a neighbouring vertex labelled \textit{Considered}.
\end{description}
\begin{definition} \label{globalaniso}
Let $\Gamma = \frac{G_{max}}{G_{min}}$ denote the \textbf{global anisotropy coefficient} where $G_{min}$ and $G_{max}$ are described in (\ref{gbounds}).
\end{definition}
\begin{description}
\item \textit{\textbf{Near Front} of $\textbf{x}_i$} ($\textbf{NF}(\textbf{x}_i)$) - Let $\textbf{x}_i$ be labelled \textit{Considered}. Define
\end{description}
\begin{equation}
\textbf{NF}(\textbf{x}_i) = \left\{\textbf{s} \in \textbf{AF} \Big| \text{ } \exists \ \widetilde{\textbf{x}} \in \textbf{s} \Big |    \left\|\widetilde{\textbf{x}} - \textbf{x}_i \right\| \leq \Gamma h_{max}\right\}.
\label{nearfront}
\end{equation}
See Figure \ref{classify}. The sets $\textbf{AF}$, $\textbf{NF}(\textbf{x}_i) \subset X_{n-1}$ change throughout the execution of the OUM due to the vertices of $X$ being relabelled from \textit{Far} to \textit{Considered} to \textit{Accepted}.

Define the discrete set of controls $\widetilde{\mathcal{U}}$ 
\begin{equation} \label{controldiscrete} \displaystyle \widetilde{\mathcal{U}} = \left\{ \widetilde{\textbf{u}}(\cdot) \in \mathcal{U} \Big | \widetilde{\textbf{u}}(t) = \widetilde{\textbf{u}}_i, \widetilde{\textbf{u}}_i \in \mathbb{S}^{n-1} \text{ while } \textbf{y}(t) \in \textbf{s} \in X \right\}. \end{equation}

The distance between vertex $\textbf{x}_i$ and $\textbf{x} \in \textbf{s} \in X_{n-1}$, where $\textbf{x} = \sum_{j=0}^{n-1} \zeta_j\textbf{x}_j^{\textbf{s}}$ is denoted $\tau_{\textbf{s}}(\textbf{x}_i,\zeta)=\left\| \sum_{j=0}^{n-1} \zeta_j\textbf{x}_j^{\textbf{s}} - \textbf{x}_i \right\| = \left\|\textbf{x}-\textbf{x}_i\right\|$. The direction from $\textbf{x}_i$ to $\textbf{x}$ is $\textbf{u}_{\textbf{s}}(\textbf{x}_i,\zeta)=\frac{\textbf{x}-\textbf{x}_i}{\tau_{\textbf{s}}(\textbf{x}_i,\zeta)}$. The update for $\textbf{x}_i$ provided by $\textbf{s} = \textbf{x}_0^{\textbf{s}}\textbf{x}_1^{\textbf{s}}\cdots\textbf{x}_{n-1}^{\textbf{s}}$ is a first-order approximation of the DPP (\ref{bellmancont}),
\begin{equation} \displaystyle
\widetilde{C}_{\textbf{s}}(\textbf{x}_i) \\ 
= \min_{\zeta \in \Xi_{n-1}} \left\{ \sum_{j = 0}^{n-1} \zeta_j \widetilde{V}(\textbf{x}_j^{\textbf{s}}) + \tau_{\textbf{s}}(\textbf{x}_i, \zeta)g(\textbf{x}_i,\textbf{u}_{\textbf{s}}(\textbf{x}_i,\zeta)) \right\},
\label{oumminimize}
\end{equation}
where $\zeta = (\zeta_0,\zeta_1,...,\zeta_{n-1}) \in \Xi_{n-1}$. The optimizing direction is captured by updating $\textbf{x}_i$ from its \textit{Near Front} \cite{oum}. The update formula over all of $\textbf{NF}(\textbf{x}_i)$ is

\begin{equation} \displaystyle
\widetilde{C}(\textbf{x}) \\
= \min_{\textbf{s} \in \textbf{NF}(\textbf{x}_i)} \widetilde{C}_{\textbf{s}}(\textbf{x}_i).
\label{oumupdate}
\end{equation}
\noindent Note that the minimizing update along all of $\textbf{NF}(\textbf{x}_i)$ (\ref{oumupdate}) does not necessarily come from $\textbf{s} \in X_{n-1}$ where $\textbf{x}_i$ is a neighbour of $\textbf{s}$. 

The algorithm can now be stated. Recall that any vertex $\textbf{x}_i \in X_0$ is labelled only one of \textit{Accepted}, \textit{Considered} or \textit{Far} at any instant of the algorithm. 
\begin{enumerate}
\item Label all vertices $\textbf{x}_i \in X_0$ \textit{Far}, assigning $\widetilde{V}(\textbf{x}_i) = K$ (where $K$ is large).
\item For each vertex $\textbf{x}_i \in X_0 \cap \Omega^c$, relabel $\textbf{x}_i$ \textit{Accepted}, and set $\widetilde{V}(\textbf{x}_i) = q(\hat{\textbf{x}})$ where $\hat{\textbf{x}} = \argmin_{\widetilde{\textbf{x}} \in \partial \Omega}\left\|\textbf{x}_i-\widetilde{\textbf{x}}\right\|$. 
\item Relabel all neighbours of \textit{Accepted} vertices $\textbf{x}_i$ that have \textit{Far} label, to \textit{Considered}. For these vertices, compute $\widetilde{V}(\textbf{x}_i) = \widetilde{C}(\textbf{x}_i)$ according to (\ref{oumupdate}).
\item Relabel vertex $\overline{\textbf{x}}_i$ with \textit{Considered} label with lowest value $\widetilde{V}(\overline{\textbf{x}}_i)$ with $\textit{Accepted}$ label. If all vertices in $X$ are labelled \textit{Accepted}, terminate the algorithm.
\item Relabel all neighbouring vertices $\textbf{x}_i$ of $\overline{\textbf{x}}_i$ with \textit{Far} label to \textit{Considered}.  For these vertices, compute $\widetilde{C}(\textbf{x}_i)$ using (\ref{oumupdate}) and set $\widetilde{V}(\textbf{x}_i) = \widetilde{C}(\textbf{x}_i)$.
\item Recompute $\widetilde{C}(\textbf{x}_i)$ for all other $\textbf{x}_i$ with \textit{Considered} label using (\ref{oumupdate}) such that $\overline{\textbf{x}}_i \in \textbf{NF}(\textbf{x}_i)$, using only $\textbf{s} \in \textbf{NF}(\textbf{x}_i)$ such that $\overline{\textbf{x}_i} \in \textbf{s}$. If $\widetilde{V}(\textbf{x}_i) > \widetilde{C}(\textbf{x}_i)$, then update $\widetilde{V}(\textbf{x}_i) = \widetilde{C}(\textbf{x}_i)$. Go to Step 4.
\end{enumerate}

The domain of $\widetilde{V}$ will be extended from $X_0$ to all of $X$. Define
\begin{equation} \label{omegax} \overline{\Omega}_X = \left\{ \bigcup_{\textbf{s} \in X_n} \bigcup_{\zeta \in \Xi_{n}} \sum_{j=0}^n \zeta_j\textbf{x}_j^{\textbf{s}}  \right\}. \end{equation} 
From (\textbf{M2}), $\overline{\Omega} \subseteq \overline{\Omega}_X$. 

The domain of the spatial dimension of value function $V$ and $g$ (and as a result $H$) are extended from $\overline{\Omega}$ to $\overline{\Omega}_X$. For $\textbf{x} \in \overline{\Omega}^c \cap \overline{\Omega}_X$, let $$\hat{\textbf{x}} = \argmin_{\widetilde{\textbf{x}}\in\partial\Omega} \left\|\textbf{x} - \widetilde{\textbf{x}}\right\|, V(\textbf{x}) = q(\hat{\textbf{x}}), \text{ and } g(\textbf{x},\textbf{u}) = g(\hat{\textbf{x}},\textbf{u}).$$

The domain of $\widetilde{V}$ is extended from $X_0$ to $\overline{\Omega}_X$ by linear interpolation using barycentric coordinates. For $\textbf{x} \in \textbf{s} = \textbf{x}_0^{\textbf{s}}\textbf{x}_1^{\textbf{s}}\cdots\textbf{x}_n^{\textbf{s}} \in X_n$, $$\displaystyle \widetilde{V}(\textbf{x}) = \sum_{j=0}^n \zeta_j \widetilde{V}(\textbf{x}_j^{\textbf{s}}), \text{ where } \textbf{x} = \sum_{j=0}^n \zeta_j\textbf{x}_j^{\textbf{s}}.$$

Most of the effort in the implementation of the OUM occurs in the maintenance and the searching of \textbf{AF} and \textbf{NF}$(\textbf{x}_i)$. The focus of this paper however is on the accuracy and its convergence to the true solution in relation to discretization properties. Additional discussion on the implementation and computational complexity of OUM can be found in \cite{oum}.

\section{Properties of the Approximated Value Function and Numerical Hamiltonian} \label{secoumprop}

An approximation of the Hamiltonian $H$  (\ref{hamdef}) known as the numerical Hamiltonian will be defined on the vertices $X_0$ of $X$. A similar numerical Hamiltonian was proposed in \cite{kaltonpaper}. As in \cite{kaltonpaper}, the numerical Hamiltonian will be shown to be both monotonic and consistent with the Hamiltonian (\ref{hamdef}). The consistency statement here resembles that in \cite{monneau}, which was given as an assumption for the half-order convergence proof for FMM. The proof of consistency relies on directional completeness introduced in \cite{kaltonpaper}.

Consider the OUM algorithm at the instant the vertex $\textbf{x}_i \in X_0 \cap \Omega$ is about to be relabelled $\textit{Accepted}$. The \textit{Near Front} of $\textbf{x}_i$ at this instant is denoted $\overline{\textbf{NF}}(\textbf{x}_i)$.
\begin{definition}
The \textbf{approximated characteristic direction} $\widetilde{\textbf{u}}_{\widetilde{\textbf{s}}}^*:  X_0 \cap \Omega \times \Xi_{n-1} \rightarrow \mathbb{S}^{n-1}$ at $\textbf{x}_i \in X_0 \cap \Omega$ from the OUM algorithm is 
\begin{equation} \widetilde{\textbf{u}}_{\widetilde{\textbf{s}}}^*(\textbf{x}_i,\widetilde{\zeta}^*) = \frac{\widetilde{\textbf{x}}^*-\textbf{x}_i}{\left\|\widetilde{\textbf{x}}^*-\textbf{x}_i\right\|}=\frac{\sum_{j=0}^{n-1}\widetilde{\zeta}_j^*\textbf{x}_j^{\widetilde{\textbf{s}}^*}-\textbf{x}_i}{\tau_{\widetilde{\textbf{s}}^*}(\textbf{x}_i,\widetilde{\zeta}_j^*)} \text{ where } \widetilde{\textbf{x}}=\sum_{j=0}^{n-1}\widetilde{\zeta}_j^*\textbf{x}_j^{\widetilde{\textbf{s}}^*} \label{oumvalueupdate} \end{equation}
where $\widetilde{\textbf{s}}^* \in \overline{\textbf{NF}}(\textbf{x}_i)$ and $\widetilde{\zeta}^* \in \Xi_{n-1}$ are the minimizers of (\ref{oumminimize}), (\ref{oumupdate}) when $\textbf{x}_i$ is labelled \textit{Accepted}.
\end{definition}
\begin{definition}
Let $\phi: X_0 \cap\Omega \rightarrow \mathbb{R}$. The \textbf{numerical Hamiltonian} $\widetilde{H}:X_0 \cap \Omega \times \mathbb{R} \rightarrow \mathbb{R}$ is  \small 
\begin{equation} \displaystyle \label{numham}
\widetilde{H}[\mathcal{S},\phi[\mathcal{S}]](\textbf{x}_i,\mu)
= - \min_{\textbf{s} \in \mathcal{S}} \min_{\zeta \in \Xi_{n-1}} \left\{ \frac{\sum_{j=0}^{n-1} \zeta_j \phi(\textbf{x}_j^{\textbf{s}}) - \mu}{\tau_{\textbf{s}}(\textbf{x}_i,\zeta)} + g(\textbf{x}_i,\textbf{u}_{\textbf{s}}(\textbf{x}_i,\zeta)) \right\},
\end{equation} \normalsize
where $\mathcal{S}\subset X_{n-1}$. 
\end{definition}
The argument $\phi[\mathcal{S}]$ of $\widetilde{H}$ denotes the use of the values of $\phi$ on the vertices that make up the $(n-1)$-simplices of $\phi[\mathcal{S}]$ in the optimization of (\ref{numham}). For notational brevity, the argument of $\phi$ will be dropped.

The numerical HJB equation for the OUM algorithm for all $\textbf{x}_i \in X_0 \cap \Omega$ is
\begin{equation} \label{numhamstatichjb} \widetilde{H}[\overline{\textbf{NF}}(\textbf{x}_i), \widetilde{V}](\textbf{x}_i,\widetilde{V}(\textbf{x}_i)) = 0.\end{equation} 

\begin{theorem} \label{altontheorem} \cite[Prop 5.3]{kaltonthesis}
Let $\mathcal{S} \subset X_{n-1}$. The solution $\mu$ to  $\widetilde{H}[\mathcal{S},\widetilde{V}](\textbf{x}_i,\mu) = 0$ with $\widetilde{H}$ defined by (\ref{numham}) is unique, and is given by 
\begin{equation} \label{equivvalue} \displaystyle \widetilde{\mu} = \min_{\textbf{s} \in \mathcal{S}} \min_{\zeta \in \Xi_{n-1}} \left\{\sum_{j=0}^{n-1} \zeta_j \widetilde{V}(\textbf{x}_j^\textbf{s}) + \tau(\textbf{x}_i, \zeta)g(\textbf{x}_i,\textbf{u}_{\textbf{s}}(\textbf{x}_i,\zeta))\right\}. \end{equation}
Furthermore, if $\widetilde{\textbf{s}}^* \in \mathcal{S}$ and $\widetilde{\zeta}^* \in \Xi_{n-1}$ are the minimizers in (\ref{numham}), then $\widetilde{\textbf{s}}^*$ and $\widetilde{\zeta}^*$ also minimize (\ref{equivvalue}).
\end{theorem}
\noindent From Theorem \ref{altontheorem}, finding the solution $\widetilde{V}(\textbf{x}_i)$ to (\ref{numhamstatichjb}) is equivalent to solving the update (\ref{oumupdate}) in the OUM algorithm for $\mathcal{S}=\overline{\textbf{NF}}(\textbf{x}_i)$. 

\begin{definition}\cite[Section 2.2]{kaltonpaper}
The set $\mathcal{S} \subseteq X_{n-1}$ is \textbf{directionally complete} for a vertex $\textbf{x}_i \in X_0$ if for all $\textbf{u} \in \mathbb{S}^{n-1}$ there exists $\textbf{x} \in \textbf{s}$ where $\textbf{s} \in \mathcal{S}$ such that 
\begin{equation*} \textbf{u} = \frac{\textbf{x} - \textbf{x}_i}{\left\|\textbf{x}-\textbf{x}_i\right\|}. \end{equation*}
\end{definition}

A subset $A \subset \mathbb{R}^n$ has no holes if its complement $A^c$ is connected.

\begin{lemma} \label{afclosedsurf}
Prior to each instance of Step 4 of the OUM algorithm, $(n-1)$-simplices of \textbf{AF} form the boundaries $\textbf{AF}_k$ of $j$ ($1 \leq k \leq j<\infty$) bounded open subsets $\Omega_{\textbf{AF}_j} \subset \overline{\Omega}_X$, such that each $\Omega_{\textbf{AF}_j}^c$ is connected and $\bigcup_{k=1}^j \textbf{AF}_k = \textbf{AF}$.  

Furthermore, if $\textbf{x}_m \in X_0 \cap\Omega_{\textbf{AF}_k}$, then 
\begin{enumerate}
\item the set of $(n-1)$-simplices $\textbf{AF}_k$ is directionally complete for $\textbf{x}_m$, and
\item $\textbf{x}_m$ is not labelled \textit{Accepted}.
\end{enumerate}
\end{lemma}

\begin{figure}
\begin{center}
\subfloat[][After relabelling $\textbf{x}_i$ \textit{Accepted}, $\textbf{AF}_3$ is no longer part of $\textbf{AF}$.]{\includegraphics[width=1.5in]{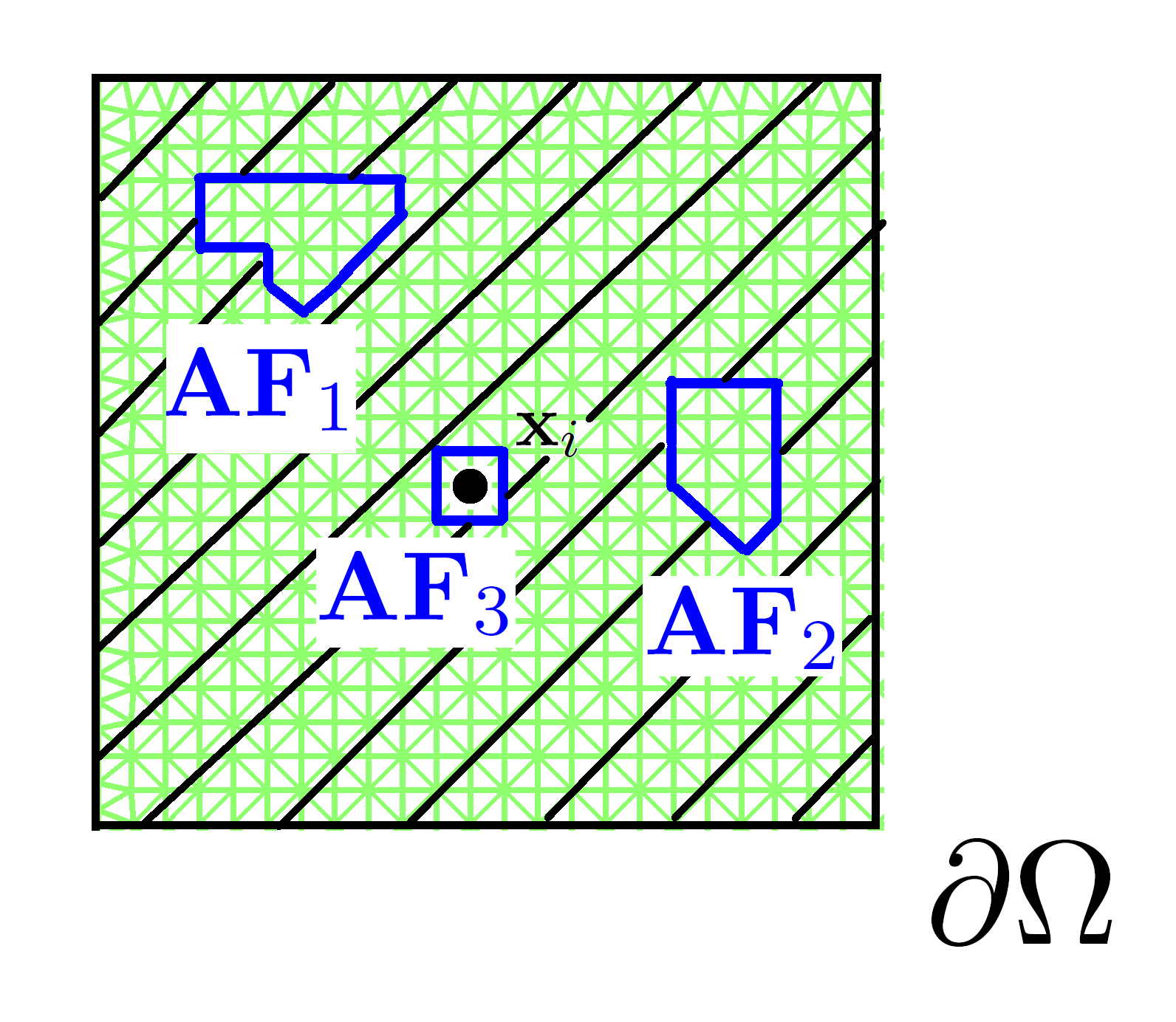}}  
\hspace{2mm} \subfloat[][After relabelling $\textbf{x}_i$ \textit{Accepted}, the other vertices in the interior of $\textbf{AF}_1$ are still not yet \textit{Accepted}.]{\includegraphics[width=1.5in]{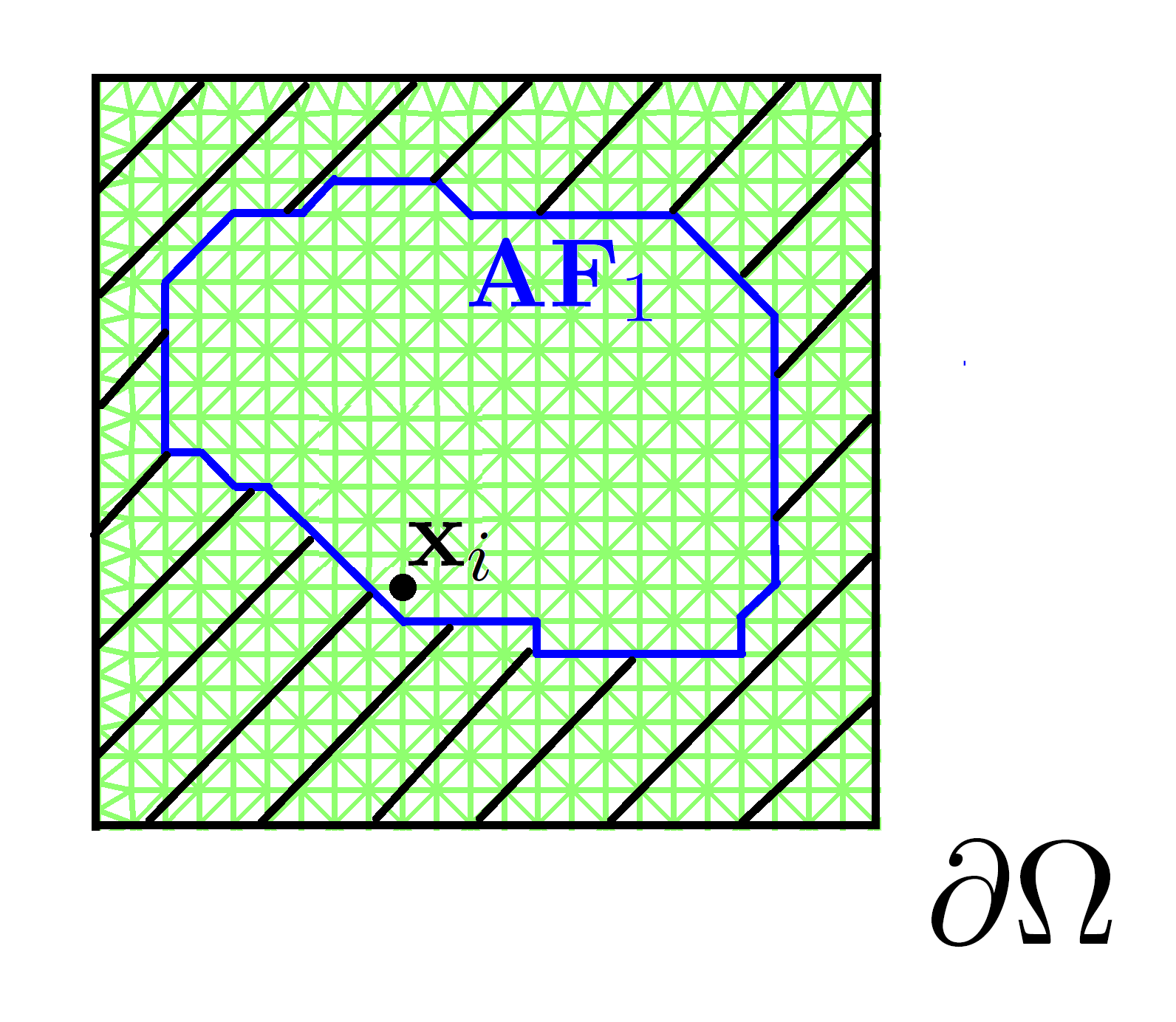}} \hspace{2mm} \subfloat[][Relabelling $\textbf{x}_i$ \textit{Accepted} splits $\textbf{AF}_1$ into two regions, each only containing not yet \textit{Accepted} vertices in their interiors.]{\includegraphics[width=1.5in]{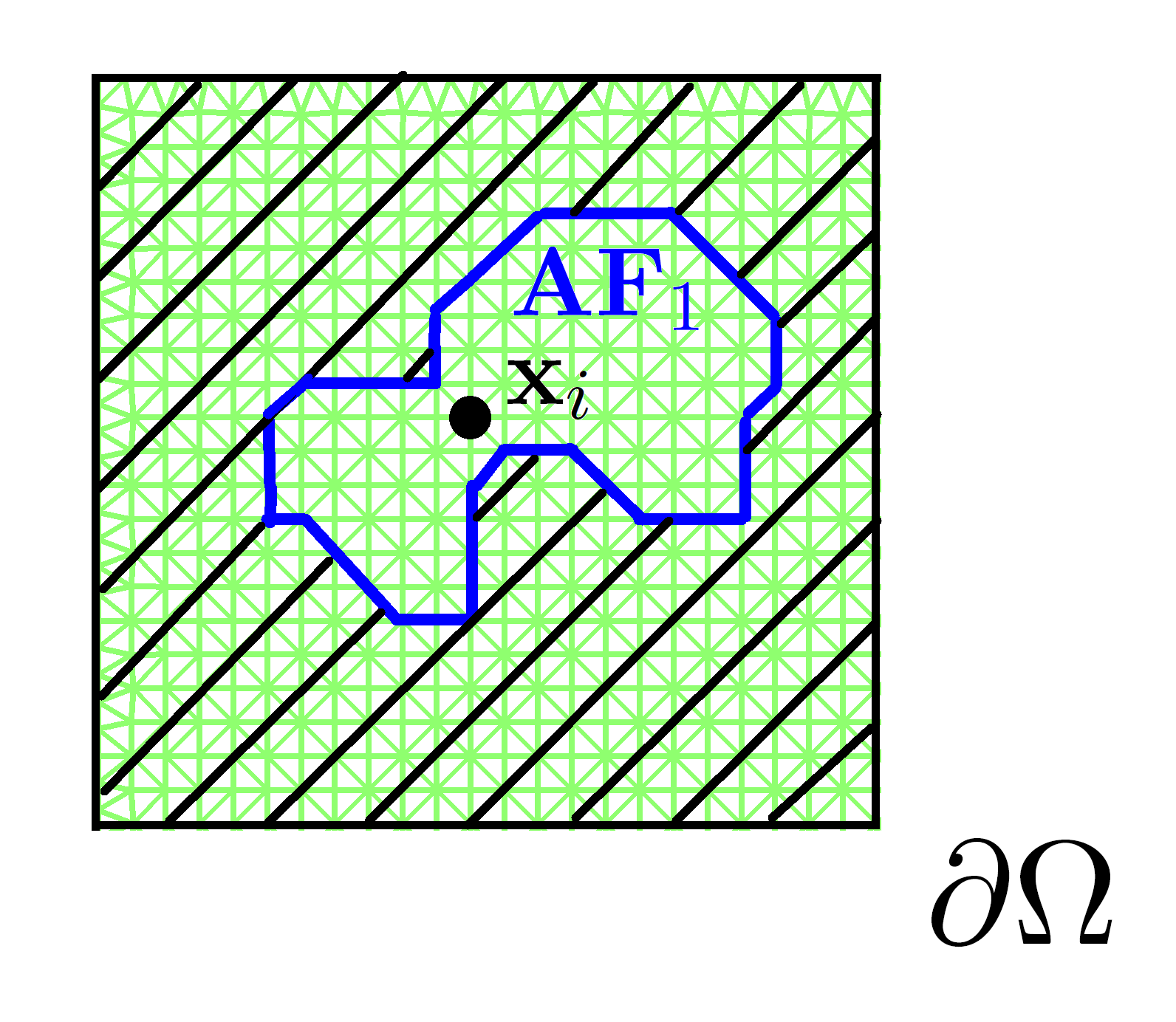}}
\end{center}
\caption{Lemma \ref{afclosedsurf}: Three cases in $\mathbb{R}^2$}
\label{afcollect}
\end{figure}

\noindent \textit{Proof.}
At the initialization (Steps 1-3) of the OUM algorithm, only vertices in $X_0 \cap\Omega^c$ are labelled \textit{Accepted}. From (\textbf{M2}) and  (\textbf{P4}), $j=1$ and $\textbf{AF}_1 = \textbf{AF}$ form a single boundary that encloses $\Omega_{\textbf{AF}_1} \supseteq \Omega $. The lemma is satisfied in the first instance of Step 4. 

The \textit{Accepted Front} and \textbf{AF} change only in Step 4 of the OUM. Proof by induction will be used. The lemma is assumed to hold prior to step 4 of the OUM. Let $\textbf{x}_i \in  X_0 \cap\Omega_{\textbf{AF}_k}$ be the vertex to be relabelled \textit{Accepted} for some $1 \leq k \leq j$. Only $\textbf{AF}_k$ and $\Omega_{\textbf{AF}_k}$ may change while $\Omega_{\textbf{AF}_{j\neq k}}$ will remain unchanged. 

If $\textbf{x}_i$ has no neighbours in $X_0 \cap\Omega_{\textbf{AF}_k}$, then the resulting $\Omega_{\textbf{AF}_k}$ and $X_0\cap\Omega_{\textbf{AF}_k}$ are both empty. See Figure \ref{afcollect}a.

If $\textbf{x}_i$ has a neighbour in $X_0 \cap\Omega_{\textbf{AF}_k}$, then $\textbf{x}_i$ is added to the \textit{Accepted Front}. If $\Omega_{\textbf{AF}_k}$ remains a single open connected subset of $\mathbb{R}^n$, $\textbf{x}_m \in X_0\cap\Omega_{\textbf{AF}_k}\backslash \{\textbf{x}_i\}$, $\textbf{AF}_k$ remains directionally complete and $\textbf{x}_m$ is not labelled \textit{Accepted}.  See Figure \ref{afcollect}b. 

Otherwise, $\Omega_{\textbf{AF}_k}$ is no longer a single open connected subset of $\mathbb{R}^n$. Thus, $\Omega_{\textbf{AF}_k}$ has been split into $p \geq 2$ non-intersecting open connected regions $\Omega_{\textbf{AF}_{k1}}$,$\Omega_{\textbf{AF}_{k2}}$,$...$,$\Omega_{\textbf{AF}_{kp}}$ with a subset of the resultant $\textbf{AF}_k$ as the boundary of each. Vertices $\textbf{x}_m \in X_0\cap\Omega_{\textbf{AF}_k}\backslash \{\textbf{x}_i\}$ are still not labelled \textit{Accepted}, and $\textbf{AF}_{kl}$ is directionally complete for $\textbf{x}_m \in \Omega_{\textbf{AF}_{kl}}$. See Figure \ref{afcollect}c. $\square$

\begin{definition} \label{dcstencil}
For every $\textbf{x}_i \in X_0\cap\Omega$, let $S(\textbf{x}_i) \subset X_{n-1}$ such that
\begin{enumerate}
\item $\overline{\textbf{NF}}(\textbf{x}_i) \subseteq S(\textbf{x}_i)$, 
\item $S(\textbf{x}_i)$ is directionally complete for $\textbf{x}_i$.
\item For all $\textbf{s} \in S(\textbf{x}_i)$, if a point $\textbf{x} \in \textbf{s}$, then $$\left\|\textbf{x} - \textbf{x}_i\right\| \leq (2\Gamma + 1)h_{max}.$$
\item $\widetilde{H}[S(\textbf{x}_i),\widetilde{V}](\textbf{x}_i,\widetilde{V}(\textbf{x}_i)) = \widetilde{H}[\overline{\textbf{NF}}(\textbf{x}_i),\widetilde{V}](\textbf{x}_i,\widetilde{V}(\textbf{x}_i))$
\end{enumerate}
\end{definition}

Such $S(\textbf{x}_i)$ will now be constructed for all $\textbf{x}_i \in X_0 \cap \Omega$ and shown to satisfy Definition \ref{dcstencil}. Let $\overline{B}_r(\textbf{x}) = \{\widetilde{\textbf{x}} \in \mathbb{R}^n | \left\|\textbf{x} - \widetilde{\textbf{x}}\right\| \leq r, r \in \mathbb{R}_+ \}$.  

\begin{definition}
Assume the OUM algorithm is at the instant that vertex $\textbf{x}_i$ labelled \textit{Considered} is about to be relabelled \textit{Accepted}.
Let $\overline{\textbf{AF}}(\textbf{x}_i)$ be the subset of $\textbf{AF}$ described in Lemma \ref{afclosedsurf} for $\textbf{x}_i$ labelled \textit{Considered}.
\end{definition}

\begin{figure}
\begin{center} \includegraphics[width=2in]{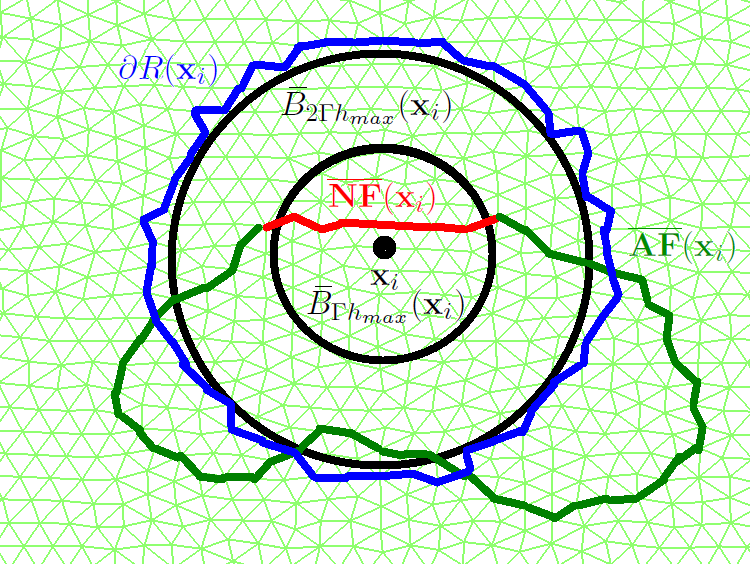} \quad \includegraphics[width=2in]{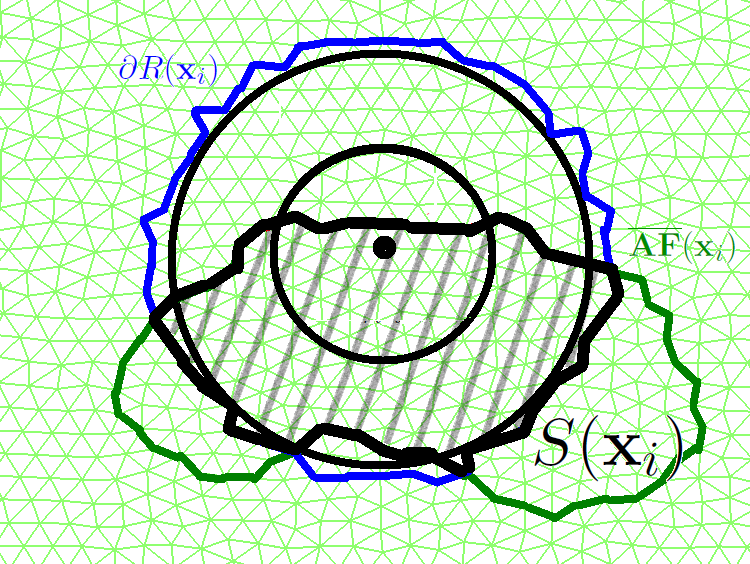}\end{center}
\caption{$S(\textbf{x}_i)$ in $\mathbb{R}^2$ - Left: Edges of $\overline{\textbf{NF}}(\textbf{x}_i)$, $\overline{\textbf{AF}}(\textbf{x}_i)$ and $\partial R(\textbf{x}_i)$ are shown. Right: $S(\textbf{x}_i)$ is the union of $\overline{\textbf{NF}}(\textbf{x}_i)$ with the boundary of the intersection of regions $R(\textbf{x}_i)$  with $\overline{\textbf{AF}}(\textbf{x}_i)$. Vertices strictly inside $S(\textbf{x}_i)$ are not labelled \textit{Accepted}.} 
\label{newNFdivide}
\end{figure}

Two cases are considered.
\begin{description}
\item \textbf{Case 1}: The set $\overline{\textbf{AF}}(\textbf{x}_i)$ lies in the interior of $\overline{B}_{2\Gamma h_{max}}(\textbf{x}_i)$, where $h_{max}$ and $\Gamma$ have been defined in Definitions \ref{hmaxdef} and \ref{globalaniso} respectively. Let $$S(\textbf{x}_i) = \overline{\textbf{AF}}(\textbf{x}_i)\cup\overline{\textbf{NF}}(\textbf{x}_i).$$ 
\item \textbf{Case 2}: Otherwise, let $R(\textbf{x}_i)$ be the region described by the smallest subset of $X_n$ in which $\overline{\Omega}_X \cup \overline{B}_{2\Gamma h_{max}}(\textbf{x}_i)$ is contained, and $\partial R(\textbf{x}_i)$ its boundary. 

%
Let $S_{\overline{\textbf{AF}}R}(\textbf{x}_i) \subset X_{n-1}$ form the boundary of the compact region $\overline{\Omega}_{\overline{\textbf{AF}}(\textbf{x}_i)} \cap R(\textbf{x}_i)$. Finally for Case 2,
\begin{equation} \label{case2} S(\textbf{x}_i) = S_{\overline{\textbf{AF}}R}(\textbf{x}_i) \cup \overline{\textbf{NF}}(\textbf{x}_i). \end{equation}
\noindent See Figure \ref{newNFdivide}. 
\end{description}

In both cases, the union with $\overline{\textbf{NF}}(\textbf{x}_i)$ ensures that $\textbf{s} \in \overline{\textbf{NF}}(\textbf{x}_i) \backslash \overline{\textbf{AF}}(\textbf{x}_i)$ are still included in $S(\textbf{x}_i)$, just as in OUM.

By construction, $S(\textbf{x}_i)$ satisfies the first three properties of Definition \ref{dcstencil}. It remains to show Property 4 in Definition \ref{dcstencil} is satisfied. 

For $\textbf{x}_i \in X_0\cap \Omega$, let $\widetilde{V}_{min}^{\textbf{AF}_{\textbf{x}_i}}$ be the minimum value on the \textit{Accepted Front} \textbf{AF} just before $\textbf{x}_i$ is labelled \textit{Accepted}. 
\begin{lemma}\cite[Lemma 7.3(i) and (iii)]{oum} \label{oumAFprop}
Assume the vertex $\textbf{x}_i \in X_0$ is about to be labelled \textit{Accepted}. Then
\begin{enumerate}
\item $\widetilde{V}_{min}^{\textbf{AF}_{\textbf{x}_i}} + h_{min}G_{min} \leq \widetilde{V}(\textbf{x}_i) \leq \widetilde{V}_{min}^{\textbf{AF}_{\textbf{x}_i}} + h_{max}G_{max}.$ 
\item If $\textbf{x}_i$ is labelled \textit{Accepted} before $\textbf{x}_j$ then $\widetilde{V}_{min}^{\textbf{AF}_{\textbf{x}_i}} \leq \widetilde{V}_{min}^{\textbf{AF}_{\textbf{x}_j}}$.
\end{enumerate}
\end{lemma}

\begin{lemma} \label{consoutside}
Let $\widetilde{\textbf{x}} = \sum_{j=0}^{n-1} \zeta_j\textbf{x}_j^{\textbf{s}}$ where $\textbf{s} = \textbf{x}_0^{\textbf{s}}\textbf{x}_1^{\textbf{s}}\cdots\textbf{x}_{n-1}^{\textbf{s}} \in X_{n-1}$, $\zeta \in \Xi_{n-1}$. If $\textbf{x}_i \in X_0$ is labelled \textit{Accepted} before all of $\textbf{x}_0^{\textbf{s}}$, $\textbf{x}_1^{\textbf{s}}$, ...,$\textbf{x}_{n-1}^{\textbf{s}}$ and $\left\|\widetilde{\textbf{x}} - \textbf{x}_i\right\| > \Gamma h_{max}$, then
\begin{equation}\label{inequality}\displaystyle  \widetilde{V}(\textbf{x}_i) < \sum_{j=0}^{n-1} \zeta_j\widetilde{V}(\textbf{x}_j^{\textbf{s}})  + \left\|\widetilde{\textbf{x}} - \textbf{x}_i\right\|g\left(\textbf{x}_i,\frac{\widetilde{\textbf{x}} - \textbf{x}_i}{\left\|\widetilde{\textbf{x}} - \textbf{x}_i\right\|}\right). \end{equation}

\end{lemma}
\noindent \textit{Proof}. From Lemma \ref{oumAFprop}, \textbf{(P2)}, Definition \ref{globalaniso} and $\widetilde{V}_{min}^{\textbf{AF}_{\textbf{x}_j^{\textbf{s}}}} < \widetilde{V}(\textbf{x}_j^{\textbf{s}})$ for $j = 1,...,n-1$,
\begin{align*} \displaystyle \widetilde{V}(\textbf{x}_i) & \leq \widetilde{V}_{min}^{\textbf{AF}_{\textbf{x}_i}} + h_{max}G_{max}, \\
 & \leq \sum_{j=0}^{n-1}\zeta_j \min\{\widetilde {V}_{min}^{\textbf{AF}_{\textbf{x}_0^{\textbf{s}}}}, \widetilde{V}_{min}^{\textbf{AF}_{\textbf{x}_1^{\textbf{s}}}},..., \widetilde{V}_{min}^{\textbf{AF}_{\textbf{x}_{n-1}^{\textbf{s}}}} \} + \Gamma h_{max} G_{min}\\
&  < \sum_{j=0}^{n-1} \zeta_j  \widetilde{V}(\textbf{x}_j^{\textbf{s}}) + \left\|\widetilde{\textbf{x}}- \textbf{x}_i\right\|g\left(\textbf{x}_i,\frac{\widetilde{\textbf{x}} - \textbf{x}_i}{\left\|\widetilde{\textbf{x}} - \textbf{x}_i\right\|}\right).
 \square
\end{align*}
\begin{lemma} \label{afsubnf} \cite[Lemma 7.1]{oum} 
Let $\textbf{x}_i$ be the vertex with \textit{Considered} label that is about to be relabelled \textit{Accepted}. Let 
\begin{equation}\label{afupdate} \displaystyle \widetilde{W}(\textbf{x}_i) = \min_{s \in \textbf{AF}} \min_{\zeta \in \Xi_{n-1}} \left\{\sum_{j = 0}^{n-1} \zeta_j\widetilde{V}(\textbf{x}_j^{\textbf{s}}) + \tau_s(\textbf{x}_i,\zeta)g(\textbf{x}_i,\textbf{u}_{\textbf{s}}(\textbf{x}_i,\zeta))\right\}. \end{equation}
Then $\widetilde{W}(\textbf{x}_i) = \widetilde{V}(\textbf{x}_i)$. 
\end{lemma}
The minimizing update from \textbf{AF} must come from $\overline{\textbf{NF}}(\textbf{x}_i)$. The next theorem states that the minimizing update $\widetilde{V}(\textbf{x}_i)$ from $S(\textbf{x}_i)$ must come from $\overline{\textbf{NF}}(\textbf{x}_i)$.
\begin{theorem} \label{NFSsame}
Let $\widetilde{V}: X_0 \rightarrow \mathbb{R}$ be computed by the OUM on mesh $X$, with weight function $g$ and boundary function $q$. Then for $\textbf{x}_i \in X_0\cap\Omega$, 
\begin{equation} \displaystyle \label{stencilupdate} \widetilde{V}(\textbf{x}_i) = \min_{\textbf{s} \in S(\textbf{x}_i)} \min_{\zeta \in \Xi_{n-1}} \left\{\sum_{j=0}^{n-1} \zeta_j \widetilde{V}(\textbf{x}_j^{\textbf{s}}) + \tau(\textbf{x}_i,\zeta)g(\textbf{x}_i, \textbf{u}_{\textbf{s}}(\textbf{x}_i,\zeta))\right\}. \end{equation}
\end{theorem}

\noindent \textit{Proof.}
Let the OUM algorithm be at the instant where vertex $\textbf{x}_i$ with \textit{Considered} label is about to be relabeled \textit{Accepted}.

Recall Case 1, where $\overline{\textbf{AF}}(\textbf{x}_i)$ is entirely inside $\overline{B}_{2\Gamma h_{max}}(\textbf{x}_i)$ and $S(\textbf{x}_i) = \overline{\textbf{AF}}(\textbf{x}_i) \cup \overline{\textbf{NF}}(\textbf{x}_i)$. Since $\overline{\textbf{AF}}(\textbf{x}_i) \subseteq \textbf{AF}$ and $\overline{\textbf{NF}}(\textbf{x}_i) \subseteq \textbf{AF}$, $S(\textbf{x}_i) \subseteq \textbf{AF}$. By Lemma \ref{afsubnf}, $\overline{\textbf{NF}}(\textbf{x}_i)$ must contain the minimizers $\widetilde{\textbf{s}}^*$ and $\widetilde{\zeta}^*$ of (\ref{stencilupdate}).

Recall Case 2, where $S(\textbf{x}_i) = S_{\overline{\textbf{AF}}R}(\textbf{x}_i) \cup \overline{\textbf{NF}}(\textbf{x}_i)$. The minimizing $\widetilde{\textbf{s}}^*$, $\widetilde{\zeta}^*$ of $S(\textbf{x}_i)$ will be shown to come from $\overline{\textbf{NF}}(\textbf{x}_i)$ by showing the updates of $S(\textbf{x}_i) \backslash \overline{\textbf{NF}}(\textbf{x}_i) = (\overline{\textbf{AF}}(\textbf{x}_i)\backslash\overline{\textbf{NF}}(\textbf{x}_i)) \cup (S(\textbf{x}_i) \cap \partial R(\textbf{x}_i))$ are at least the value from OUM. By Lemma \ref{afsubnf}, the minimizers are not from $\overline{\textbf{AF}}(\textbf{x}_i)\backslash\overline{\textbf{NF}}(\textbf{x}_i)$. 

It remains to show that updates (\ref{oumminimize}) from $\textbf{s}  \in S(\textbf{x}_i) \cap \partial R(\textbf{x}_i)$ (which are just outside $\overline{B}_{2\Gamma h_{max}}(\textbf{x}_i)$) are at least the value obtained from OUM. Because vertices of $\textbf{s}$ lie on or inside $\overline{\textbf{AF}}(\textbf{x}_i)$, they must either be on the \textit{Accepted Front} or not yet $\textit{Accepted}$ (Lemma \ref{afclosedsurf}). Three cases are considered. 
\begin{enumerate}
\item If none of the vertices of $\textbf{s}$ have been labelled \textit{Accepted}, Lemma \ref{consoutside} applies. The update for $\textbf{x}_i$ from $\textbf{s} \in S(\textbf{x}_i) \cap \partial R(\textbf{x}_i)$ is greater than $\widetilde{V}(\textbf{x}_i)$ from OUM.
\item If the vertices of $\textbf{s}$ are all on the \textit{Accepted Front}, then $\textbf{s} \in \textbf{AF}$ and Lemma \ref{afsubnf} applies. The update from $\textbf{s}$ is at least $\widetilde{V}(\textbf{x}_i)$ from OUM. 
\item If at least one but not all the vertices of $\textbf{s}$ are on the \textit{Accepted Front}, then the rest of the vertices on $\textbf{s}$ (that are not labelled \textit{Accepted}) must be labelled \textit{Considered}. Let the \textit{Accepted} and \textit{Considered} vertices of $\textbf{s}$ be denoted $\{\textbf{x}_1^\textbf{sa},...,\textbf{x}_l^\textbf{sa} \}$ and $\{\textbf{x}_1^\textbf{sc},...,\textbf{x}_k^\textbf{sc} \}$ respectively. Let $\textbf{s}$ be rewritten\\
\noindent $\textbf{s} = \textbf{x}_1^\textbf{sa}\cdots\textbf{x}_l^\textbf{sa}\textbf{x}_1^\textbf{sc}	\cdots\textbf{x}_k^\textbf{sc}$ where $l+k = n$ since $\textbf{s}$ has $n$ vertices. Let $\zeta = (\zeta_1^{\textbf{sa}},...,\zeta_l^{\textbf{sa}},\zeta_1^{\textbf{sc}},...,\zeta_k^{\textbf{sc}})$ be the barycentric coordinates for $\textbf{x} \in \textbf{s}$. By Lemma \ref{oumAFprop}, $\widetilde{V}(\textbf{x}_i) > \widetilde{V}^{\textbf{AF}_{\textbf{x}_i}}_{min}$ and Definition \ref{globalaniso}, for all $1 \leq j \leq k$,
\begin{equation*} \widetilde{V}(\textbf{x}_i) \leq \widetilde{V}^{\textbf{AF}_{\textbf{x}_i}}_{min} + h_{max}G_{max} < \widetilde{V}(\textbf{x}_j^{\textbf{sc}}) + \Gamma h_{max} G_{min}. \end{equation*}
For all $1 \leq j \leq k$, and $1 \leq m\leq l$, $\textbf{x}_j^{\textbf{sc}}$ is labelled \textit{Considered} and $\textbf{x}_m^{\textbf{sa}}$ is on its \textit{Near Front} $\overline{\textbf{NF}}(\textbf{x}_j^{\textbf{sc}})$. Thus, 
\begin{equation*}  \widetilde{V}(\textbf{x}_j^{\textbf{sc}}) \leq \widetilde{V}(\textbf{x}_m^{\textbf{sa}}) + \left\|\textbf{x}_m^{\textbf{sa}} - \textbf{x}_j^{\textbf{sc}}\right\|g\left(\textbf{x}_j^{\textbf{sc}},\frac{\textbf{x}_m^{\textbf{sa}} - \textbf{x}_j^{\textbf{sc}}}{\left\|\textbf{x}_m^{\textbf{sa}} - \textbf{x}_j^{\textbf{sc}}\right\|}\right) \leq \widetilde{V}(\textbf{x}_m^{\textbf{sa}}) + \Gamma h_{max} G_{min} \end{equation*}
\begin{equation*} \widetilde{V}(\textbf{x}_m^{\textbf{sa}}) \geq  \widetilde{V}(\textbf{x}_j^{\textbf{sc}}) - \Gamma h_{max}G_{min} > \widetilde{V}(\textbf{x}_i) - 2\Gamma h_{max}G_{min}. \end{equation*}
Consider the update for $\textbf{x}_i$ (\ref{oumminimize}) from $\textbf{s} \in S(\textbf{x}_i)\cap\partial R(\textbf{x}_i)$. For any $\zeta \in \Xi_{n-1}$,
\begin{align*}\displaystyle
& \sum_{j=0}^{n-1} \zeta_j \widetilde{V}(\textbf{x}_j^{\textbf{s}})  + \tau_{\textbf{s}}(\textbf{x}_i,\zeta) g(\textbf{x}_i,\textbf{u}_{\textbf{s}}(\textbf{x}_i,\zeta))\\
= & \left(\sum_{m=1}^l \zeta_m^{\textbf{sa}} \widetilde{V}(\textbf{x}_m^{\textbf{sa}}) \right) + \left(\sum_{j=1}^k \zeta_j^{\textbf{sc}} \widetilde{V}(\textbf{x}_j^{\textbf{sc}})\right) + \tau_{\textbf{s}}(\textbf{x}_i,\zeta)g(\textbf{x}_i,\textbf{u}_{\textbf{s}}(\textbf{x}_i,\zeta)),\\
> & \sum_{m=1}^l \zeta_m^{\textbf{sa}} (\widetilde{V}(\textbf{x}_i)- 2\Gamma h_{max}G_{min}) + \sum_{j=1}^k \zeta_j^{\textbf{sc}}(\widetilde{V}(\textbf{x}_i) - \Gamma h_{max}G_{min}) + 2\Gamma h_{max}G_{min},\\
\geq & \widetilde{V}(\textbf{x}_i),
\end{align*}
since $\tau_{\textbf{s}}(\textbf{x}_i,\zeta) > 2\Gamma h_{max}$ and ($\sum_{m=1}^l \zeta_m^{\textbf{sa}}) + (\sum_{j=1}^k \zeta_j^{\textbf{sc}}) = 1$.
\end{enumerate}
Therefore $\textbf{s} \in S(\textbf{x}_i) \cap \partial R(\textbf{x}_i)$ provides an update larger or equal to OUM. By Lemma \ref{afsubnf}, a minimizing update (\ref{stencilupdate}) in $S(\textbf{x}_i)$ must always come from $\overline{\textbf{NF}}(\textbf{x}_i)$. $\square$ \endproof

By Theorems \ref{altontheorem} and \ref{NFSsame}, $$\widetilde{H}[S(\textbf{x}_i),\widetilde{V}](\textbf{x}_i,\widetilde{V}(\textbf{x}_i)) = \widetilde{H}[\overline{\textbf{NF}}(\textbf{x}_i),\widetilde{V}](\textbf{x}_i,\widetilde{V}(\textbf{x}_i)).$$ Therefore, for $\textbf{x}_i \in  X_0 \cap \Omega$, $S(\textbf{x}_i)$ satisfies Definition \ref{dcstencil}.

The monotonicity and consistency of the numerical Hamiltonian will now be discussed.

\begin{theorem} (\textit{Monotonicity}) \cite[Proposition 2.1]{kaltonpaper} \label{monotonic}
For $\underline{\phi},\overline{\phi}:X_0 \rightarrow \mathbb{R}$ that satisfy $\underline{\phi}(\textbf{x}_j) \leq \overline{\phi}(\textbf{x}_j)$ for all $\textbf{x}_j \in  X_0 \cap \Omega$, and $\phi(\textbf{x}_i) = \overline{\phi}(\textbf{x}_i) = \underline{\phi}(\textbf{x}_i) \in \mathbb{R}$,
$$\widetilde{H}[S(\textbf{x}_i),\underline{\phi}](\textbf{x}_i,\phi(\textbf{x}_i)) \geq \widetilde{H}[S(\textbf{x}_i),\overline{\phi}](\textbf{x}_i,\phi(\textbf{x}_i)).$$
\end{theorem}

\begin{theorem} (\textit{Consistency})
There exists $C_1 \in \mathbb{R}_+$ (not dependent on $h_{max}$) for all $\textbf{x}_i \in X_0 \cap \Omega$ and  $\phi \in C^2(\Omega)$, such that
\begin{equation*}
|H(\textbf{x}_i,\nabla \phi)-\widetilde{H}[S(\textbf{x}_i),\phi](\textbf{x}_i,\phi(\textbf{x}_i))| \leq C_1\left\|\nabla^2\phi\right\|_2 h_{max}.
\end{equation*}
where $\left\| A \right\|_2$ is the maximum singular value of $A \in \mathbb{R}^{n\times n}$.
\label{consistent}
\end{theorem}

\noindent \textit{Proof.}
Let $\phi \in C^2(\Omega)$ and $\textbf{x}_i \in X_0 \cap \Omega$. Recall $S(\textbf{x}_i)$ is directionally complete, so the characteristic direction (Definition \ref{chardirdef}) $\textbf{u}^*$ can be described using barycentric coordinates $\zeta^* =(\zeta_0^*,\zeta_1^*,...,\zeta_{n-1}^*)\in \Xi_{n-1}$ from an appropriate simplex $\textbf{s}^* \in S(\textbf{x}_i)$. Let $\textbf{x}^*=\sum_{j=0}^{n-1}\zeta^*_j\textbf{x}^{\textbf{s}^*}_j$. Taylor's theorem will be used on $H$ (\ref{hamdef}). Let $\textbf{c}^*$ and $\textbf{c}_j^{*}$ for $j=0,1,...,n-1$ denote the points arising from Taylor's theorem  on the line segments between $\textbf{x}^*$ and $\textbf{x}_i$ and $\textbf{x}^*$ and $\textbf{x}^{\textbf{s}^*}_j$ respectively. Since $\sum_{j=0}^{n-1}\zeta^*_j\nabla\phi(\textbf{x}^*)^T(\textbf{x}^{\textbf{s}^*}_j-\textbf{x}^*)=\nabla\phi(\textbf{x}^*)^T(\textbf{x}^*-\textbf{x}^*)=0$, evaluating both $H$ and $\widetilde{H}$ at $\textbf{s}^*$ and $\zeta^*$,
\begin{align*}
 &H(\textbf{x}_i,\nabla \phi) - \widetilde{H}[S(\textbf{x}_i),\phi](\textbf{x}_i,\phi(\textbf{x}_i)) \\
& \leq -\frac{\sum_{j=0}^{n-1} \frac{{\zeta}_j^*}{2}(\textbf{x}_j^{{\textbf{s}}^*} - {\textbf{x}}^*)^T\nabla^2\phi(\textbf{c}_j^*)(\textbf{x}_j^{{\textbf{s}}^*} - {\textbf{x}}^*) +\frac{1}{2}({\textbf{x}}^* - \textbf{x}_i)^T\nabla^2\phi(\textbf{c}^*)({\textbf{x}}^* - \textbf{x}_i )}{\tau_{{\textbf{s}}^*}(\textbf{x}_i,{\zeta}^*)},\\
& \leq \frac{1}{h_{min}}\left(\frac{\sum_{j=0}^{n-1} \zeta_j^*}{2} \left\|\nabla^2 \phi\right\|_{2} h_{max}^2 + \frac{1}{2} \left\|\nabla^2\phi\right\|_{2} (2\Gamma+1)^2 h_{max}^2\right), \\
& \leq \frac{M}{2}\left\|\nabla^2\phi\right\|_{2}(1 + (2\Gamma+1)^2)h_{max},
\end{align*}
since the point $\textbf{x}^* \in \textbf{s}^* \in S(\textbf{x}_i)$ is at most $(2\Gamma + 1)h_{max}$ from $\textbf{x}_i$ and at most $h_{max}$ away from any of the vertices of $\textbf{s}^*$. The distance $\tau_{\textbf{s}^*}(\textbf{x}_i,\zeta^*)$ is at least the minimum simplex height $h_{min}$ and $M$ from $(\textbf{M1})$ satisfies $1 \leq \frac{h_{max}}{h_{min}} \leq M$.
The proof for $\widetilde{H}[S(\textbf{x}_i),\phi](\textbf{x}_i,\phi(\textbf{x}_i)) - H(\textbf{x}_i,\nabla \phi)$ yields the same estimate using the minimizers of $\widetilde{H}$, $(n-1)$-simplex $\widetilde{\textbf{s}}^* \in S(\textbf{x}_i)$ and $\widetilde{\zeta}^* \in \Xi_{n-1}$. The theorem is proved with $C_1 = \frac{M}{2}(1+(2\Gamma+1)^2)$. $\square$ \endproof
A similar consistency property was assumed in \cite{monneau} for the half-order proof for FMM. A similar proof without rate using similar arguments was given in \cite[Prop 2.2]{kaltonpaper} for the Monotone Acceptance OUM.
\section{OUM Error Bound} \label{secerrorbound}
The error bound proof will be presented. Several definitions and results are first required.

\begin{lemma} \cite{angellnotes}
\label{obtuseconvex}
Let $\textbf{x} \in \mathbb{R}^n$. If $\overline{\Omega}$ is convex, then $\textbf{z}^*  = \argmin_{\textbf{z} \in \overline{\Omega}} \left\|\textbf{x} - \textbf{z}\right\|$  is unique, and satisfies
\begin{equation} \label{lemma61}
( \textbf{x}-\textbf{z}^*) \cdot (\textbf{w} - \textbf{z} ^*) \leq 0 , \text{ for all } \textbf{w} \in \overline{\Omega}.
\end{equation}
\end{lemma}

\begin{lemma} \label{valuelipschitz}
The value function $V$ is globally Lipschitz-continuous over $\overline{\Omega}_X$. That is, there exists $L_V \in \mathbb{R}_+$ such that for any $\textbf{x}_1, \textbf{x}_2 \in \overline{\Omega}_X$,
\begin{equation*} |V(\textbf{x}_1) - V(\textbf{x}_2)| \leq L_V\left\|\textbf{x}_1 - \textbf{x}_2\right\|.\end{equation*} 
\end{lemma}

An outline of the proof is given using three cases.\\
\noindent \textbf{Case 1: $\textbf{x}_1, \textbf{x}_2 \in \overline{\Omega}_X \cap \Omega^c$.} This is an exercise in \cite[Exercise 2.8d]{borweinlewis}, which can be shown using the Cauchy-Schwartz inequality and Lemma \ref{qlipschitz}.

\noindent \textbf{Case 2}: $\textbf{x}_1, \textbf{x}_2 \in \Omega$. This is shown in  \cite[Lemma 2.2.7]{vladthesis} with constant $G_{max}$.

\noindent \textbf{Case 3}: $\textbf{x}_1 \in \Omega $ and $\textbf{x}_2 \in \overline{\Omega}_X \cap \Omega^c$. This can be shown using Lemma \ref{obtuseconvex} and 
\begin{equation*} L(\textbf{a},\textbf{b}) \leq L(\textbf{a},\textbf{c}) + L(\textbf{c},\textbf{b}),  \end{equation*}
for $\textbf{a}, \textbf{b},\textbf{c} \in \overline{\Omega}$. 
For $\textbf{x}_1,\textbf{x}_2 \in \overline{\Omega}_X$, a valid Lipschitz constant is $L_V=2G_{max}$.

\begin{lemma}\cite[Lemma 2.2.9]{vladthesis} \label{Vbound}
Let $\textbf{x} \in \overline{\Omega}_X$ . Let $\widetilde{\textbf{x}} = \argmin_{\textbf{z} \in \partial \Omega} |\textbf{x} - \textbf{z}|$. The value function $V$ satisfies 
\begin{equation*} q_{min} \leq V(\textbf{x}) \leq G_{max}\left\|\textbf{x} - \widetilde{\textbf{x}}\right\| + q_{max}. \end{equation*}
\end{lemma}
The proof is shown in \cite{vladthesis} for $\textbf{x} \in \overline{\Omega}$. The proof is trivial for $\textbf{x} \in \overline{\Omega}_X \cap \overline{\Omega}^c$. 

\begin{lemma}\cite[Lemma 7.5]{oum} \label{aVlipschitz}
Let $\widetilde{V}: X_0 \rightarrow \mathbb{R}$ obtained by the Ordered Upwind Method. There exists $L_{\widetilde{V}} \in \mathbb{R}_+$ for any $\textbf{x}_i, \textbf{x}_j \in X_0$,  such that
\begin{equation*} |\widetilde{V}(\textbf{x}_i) - \widetilde{V}(\textbf{x}_j)| \leq L_{\widetilde{V}} |\textbf{x}_i - \textbf{x}_j|. \end{equation*}
\end{lemma}
A possible Lipschitz constant for $\widetilde{V}$ is $L_{\widetilde{V}} = M^2 G_{max}$ \cite{oum}, where $M$ is described in (\textbf{M1}). Similar proof from case 1 and case 3 of Lemma \ref{valuelipschitz} is valid with a restriction of $\textbf{x} \in X_0$ and function $L$ (\ref{Lequation}) is replaced with $\widetilde{L}: X_0 \times X_0 \rightarrow \mathbb{R}$,
\begin{equation} \label{Lequationapprox} \displaystyle \widetilde{L}(\textbf{x}_1,\textbf{x}_2) = \inf_{\textbf{u}(\cdot) \in \widetilde{\mathcal{U}}} \left\{ \int_0^{\tau} g(\textbf{y}_{\textbf{x}_1}(s),\textbf{u}(s))ds \text{ } \Big | \text{ }\textbf{y}_{\textbf{x}_1}(\tau) = \textbf{x}_2, \textbf{y}_{\textbf{x}_1}(t) \in \overline{\Omega}, t \in (0,\tau) \right\}. \end{equation}
where $\widetilde{\mathcal{U}}$ is defined in (\ref{controldiscrete}).
\begin{lemma}\cite[Lemma 7.2]{oum} \label{aVbound}
Let $\textbf{x} \in \textbf{s}$ where $\textbf{s} \in X_n$ and $\widetilde{\textbf{x}} = \argmin_{\textbf{z} \in \partial \Omega} \left\| \textbf{x} - \textbf{z}\right\|.$ Then 
\begin{equation*}  q_{min} \leq \widetilde{V}(\textbf{x}) \leq G_{max}|\textbf{x} - \widetilde{\textbf{x}}| + q_{max}. \end{equation*}
\end{lemma}
The proof is shown in \cite{oum} for $\textbf{x} \in \overline{\Omega}$. The proof is trivial for $\textbf{x} \in \overline{\Omega}^c$. 

The next lemma states that any point on the boundary $\partial \Omega$ must be at most $h_{max}$ away from its nearest vertex of $X$ outside of $\Omega$.
\begin{lemma} \label{boundarylemma}
If $\textbf{x} \in \partial \Omega$, there exists $\textbf{x}_i \in X_0 \cap \Omega^c$ such that
\begin{equation} \left\|\textbf{x} - \textbf{x}_i\right\| \leq h_{max}. \end{equation}
\end{lemma}

\noindent \textit{Proof.} Assumption  (\textbf{M2}) states that $\overline{\Omega}$ is contained in $X$.  The point $\textbf{x} \in \textbf{s}$ where $\textbf{s} \in X_n$.  Since $\overline{\Omega}$ is convex (\textbf{P4}), and $\textbf{x}$ can be described by barycentric coordinates of $\textbf{s}$, at least one of the vertices of $\textbf{s}$ must be outside $\Omega$. Furthermore, for all $1 \leq j \leq n$,
\begin{equation*}\left\|\textbf{x} - \textbf{x}_j^s\right\| \leq \max_{1\leq k \leq n} \left\|\textbf{x}_k^{\textbf{s}} - \textbf{x}_j^{\textbf{s}}\right\| \leq h_{max}.  \square \end{equation*}

The following definitions provide a weaker description of the gradient for functions that are not necessarily differentiable. 
Let $A$ be a bounded subset of $\mathbb{R}^n$.

\begin{definition} \label{subgraddef}
The vector $\textbf{p} \in \mathbb{R}^n$ is a \textbf{subgradient} of a function $f: A \rightarrow \mathbb{R}$ at $\textbf{x}_0 \in A$ if there exists $\delta > 0$ such that for any $\textbf{x} \in B_{\delta}(\textbf{x}_0)$,
\begin{equation*} f(\textbf{x}) - f(\textbf{x}_0) \geq \textbf{p} \cdot (\textbf{x} - \textbf{x}_0). \end{equation*}
\end{definition}

\begin{definition} \label{supergraddef}
The vector $\textbf{p} \in \mathbb{R}^n$ is a \textbf{supergradient} of a function $f: A \rightarrow \mathbb{R}$ at $\textbf{x}_0 \in A$ if there exists $\delta > 0$ such that for any $\textbf{x} \in B_{\delta}(\textbf{x}_0)$,
\begin{equation*} f(\textbf{x}) - f(\textbf{x}_0) \leq \textbf{p} \cdot (\textbf{x} - \textbf{x}_0). \end{equation*}
\end{definition}
Let $D^-f(\textbf{x}_0)$ and $D^+f(\textbf{x}_0)$ denote the sets of all subgradients and supergradients of $f$ at $\textbf{x}_0$ respectively.

\begin{lemma} \label{onesidedlemma}
Let $f: A \rightarrow \mathbb{R}$ be globally Lipschitz-continuous with Lipschitz constant $C$ and $\textbf{x}_0 \in A$. If $\textbf{p} \in D^-f(\textbf{x}_0)\cup D^+f(\textbf{x}_0)$	, then $$ \left\|\textbf{p}\right\| \leq C .$$
\end{lemma}

\noindent \textit{Proof.}
Let $\textbf{x}_0 \in A$, $\textbf{b} \in \mathbb{S}^{n-1}$, $\delta > 0$, such that $\textbf{x}_0 + \delta\textbf{b} \in A$.  Let $\textbf{p} \in D^-f(\textbf{x}_0)$ (Definition \ref{subgraddef}). The Lipschitz continuity of $f$ gives
\begin{equation*} C\left\|\textbf{x}_0 + \delta\textbf{b} - \textbf{x}_0\right\| \geq f(\textbf{x}_0 + \delta\textbf{b}) - f(\textbf{x}_0) \geq \textbf{p} \cdot (\textbf{x}_0 + \delta\textbf{b} - \textbf{x}_0). \end{equation*}
Choosing $\textbf{b}=\frac{\textbf{p}}{\left\|\textbf{p}\right\|}$ gives $\left\|\textbf{p}\right\| \leq C$. The proof is analogous for $\textbf{p} \in D^+f(\textbf{x}_0)$. $\square$ 

\begin{lemma} \label{equivalentsub} \cite[Lemma 1.7]{bardi}
A vector $\textbf{p} \in D^-f(\textbf{x}_0)$ if and only if there exists $\phi \in C^1(\Omega) \rightarrow \mathbb{R}$ such that $f - \phi$ has a local minimum at $\textbf{x}_0$.
Similarly, a vector $\textbf{p} \in D^+f(\textbf{x}_0)$ if and only if there exists $\phi \in C^1(\Omega) \rightarrow \mathbb{R}$ such that $\nabla \phi(\textbf{x}_0) = \textbf{p}$, and $f - \phi$ has a local maximum at $\textbf{x}_0$.
\end{lemma}

The approximated value function $\widetilde{V}$ is in a sense a viscosity solution for the numerical HJB equation (\ref{numhamstatichjb}). 

\begin{definition}	 \label{subnumham}
Let $\hat{\textbf{x}} = \argmin_{\textbf{x} \in \partial \Omega} \left\|\textbf{x}_i - \textbf{x}\right\|$. A \textbf{subsolution of the numerical HJB equation} (\ref{numhamstatichjb}) $\underline{\widetilde{V}}: X_0 \rightarrow \mathbb{R}$  satisfies
\begin{displaymath}
   \left\{
     \begin{array}{lr}
       \underline{\widetilde{V}}(\textbf{x}_i) \leq q(\hat{\textbf{x}}) & \text{ for } \textbf{x}_i \in  X_0 \cap \Omega^c,\\
       \widetilde{H}[\overline{\textbf{NF}}(\textbf{x}_i),\underline{\widetilde{V}}](\textbf{x}_i,\underline{\widetilde{V}}(\textbf{x}_i)) \leq 0 & \text{ for } \textbf{x}_i \in  X_0 \cap \Omega.
     \end{array}
   \right.
\end{displaymath} 
\end{definition}

\begin{definition} \label{supnumham}
Let $\hat{\textbf{x}} = \argmin_{\textbf{x} \in \partial \Omega} \left\|\textbf{x}_i - \textbf{x}\right\|$. A \textbf{supersolution of the numerical HJB equation} (\ref{numhamstatichjb}) $\widetilde{\overline{V}}: X_0 \rightarrow \mathbb{R}$ satisfies 
\begin{displaymath}
   \left\{
     \begin{array}{lr}
		\widetilde{\overline{V}}(\textbf{x}_i) \geq q(\hat{\textbf{x}}) & \text{ for } \textbf{x}_i \in X_0 \cap \Omega^c,
\\
              \widetilde{H}[\overline{\textbf{NF}}(\textbf{x}_i),\widetilde{\overline{V}}](\textbf{x}_i,\widetilde{\overline{V}}(\textbf{x}_i)) \geq 0 & \text{ for } \textbf{x}_i \in  X_0 \cap \Omega.
     \end{array}
   \right.
\end{displaymath} 
\end{definition}

\begin{definition} \label{numhamdef} A \textbf{solution of the numerical HJB equation} (\ref{numhamstatichjb})  $\widetilde{V}$ is both a subsolution and a supersolution of the numerical HJB equation (\ref{numhamstatichjb}). \end{definition}

By Theorem \ref{altontheorem}, the approximate value function $\widetilde{V}$ produced by the OUM algorithm is a solution of the numerical HJB equation. Hence, it is both a subsolution and supersolution of the numerical HJB equation. Recall the definition of $\overline{\Omega}_X$ (\ref{omegax}).

\begin{theorem}
\label{oumerror}
Let $V: \overline{\Omega}_X \rightarrow \mathbb{R}$ be a viscosity solution of (\ref{conttheo1}) and $\widetilde{V}: X_0 \rightarrow \mathbb{R}$ be a solution of the numerical HJB equation (\ref{numhamstatichjb}). There exist $C, h_0 > 0$, both independent of $h_{max}$ such that
\begin{equation} \displaystyle \max_{\textbf{x}_i \in X_0} |V(\textbf{x}_i) - \widetilde{V}(\textbf{x}_i)| \leq C \sqrt{h_{max}}, 
\label{mainresult}
\end{equation} for every $\textbf{x}_i \in X_0$ and $h_{max} < h_0$.
\end{theorem}

\noindent \textit{Proof.}
The proof is trivial for $\textbf{x}_i \in X_0\cap\overline{\Omega}^c$. Otherwise, $\textbf{x}_i \in X_0 \cap \overline{\Omega}$. 
Since $\overline{\Omega} \subseteq \mathbb{R}^n$ is bounded, define 
 \begin{equation} \label{distancebound} d_{\Omega} = \max_{\textbf{x},\widetilde{\textbf{x}} \in \partial \Omega} \left\|\textbf{x} - \widetilde{\textbf{x}}\right\|, \end{equation}
 \begin{equation} \label{overallbound} C_0 = \max \{L_V, L_{\widetilde{V}}, |q_{min}|, d_{\Omega} G_{max} + |q_{max}| \},\end{equation} where $L_V, L_{\widetilde{V}}, |q_{min}|, G_{max}$, and $|q_{max}|$ are from Lemmas \ref{valuelipschitz}, \ref{Vbound}, \ref{aVlipschitz}, \ref{aVbound}. 

For $\textbf{x}_i \in X_0 \cap \overline{\Omega}$, the result of the theorem is shown for $V(\textbf{x}_i) - \widetilde{V}(\textbf{x}_i)$. A similar argument for $\widetilde{V}(\textbf{x}_i) - V(\textbf{x}_i)$ can be made.

Two parameters $\epsilon$ and $\lambda$ are used to determine the error bound. For $\epsilon > 0$ and $0 < \lambda < 1$, define $\Phi: \overline{\Omega} \times X_0 \rightarrow \mathbb{R}$
\begin{equation} \label{phidef} \Phi(\textbf{x},\textbf{x}_i) = \lambda V(\textbf{x}) - \widetilde{V}(\textbf{x}_i) - \frac{\left\|\textbf{x} - \textbf{x}_i\right\|^2}{2\epsilon}.   \end{equation}
Let $\overline{\textbf{x}}  \in \overline{\Omega}$ and $\overline{\textbf{x}}_i \in X_0$ maximize $\Phi$, over the compact set $\overline{\Omega} \times X_0$. Define
\begin{equation}\label{Mdef} \displaystyle M_{\epsilon,\lambda} = \max_{\textbf{x} \in \overline{\Omega}, \textbf{x}_i \in X_0} \Phi(\textbf{x},\textbf{x}_i) = \Phi(\overline{\textbf{x}},\overline{\textbf{x}}_i). \end{equation}

For $\textbf{x}_i \in X_0 \cap \overline{\Omega}$,  using (\ref{phidef}) and (\ref{Mdef}) with $V(\textbf{x}_i) \leq C_0$ from Lemma \ref{Vbound} (boundedness of $V$),
\begin{equation} \label{boundO} V(\textbf{x}_i) - \widetilde{V}(\textbf{x}_i) \leq (1-\lambda)  V(\textbf{x}_i) + M_{\epsilon,\lambda} \leq  C_0(1-\lambda) + M_{\epsilon,\lambda}. \end{equation} 
Choose $\lambda$ such that 
\begin{equation} \label{lambdachoice}1-\lambda = \frac{2}{G_{min}}\left(\frac{C_1}{\epsilon}h_{max} + C_0L_g\epsilon \right), \end{equation}
where $L_g$ is defined in (\ref{lipschitzweight}), and $C_1 = \frac{M(1+(2\Gamma + 1)^2)}{2}$ is defined in Theorem \ref{consistent} with $M$ in (\textbf{M1}) and  $\Gamma = \frac{G_{max}}{G_{min}}$. 

The result of the theorem will be true with $\epsilon = \sqrt{h_{max}}$. Therefore, it is sufficient to pick $h_0$ small enough so that for all $h_{max} < h_0$, $0 < (1-\lambda) < 1$ is satisfied. Setting (\ref{lambdachoice}) less than $1$, with $\epsilon = \sqrt{h_{max}}$ yields
$h_{max} < \frac{G_{min}^2}{4(C_1+C_0L_g)^2}$. Let $h_0 = \min\{\frac{G_{min}^2}{4(C_1+C_0L_g)^2},1\}$. 

The point $\overline{\textbf{x}}$ in (\ref{phidef}) must belong to $\Omega$ or $\partial \Omega$, while $\overline{\textbf{x}}_i$ must belong to $X_0\cap\Omega$ or $X_0\cap\Omega^c$.  
An outline of the remainder of  proof is as follows. 
\begin{description}
\item \textbf{Step 1}: Show that at most only one of $\overline{\textbf{x}}$ and $\overline{\textbf{x}}_i$ may be in $\Omega$. 
\item \textbf{Step 2}: Find an upper bound for $M_{\epsilon,\lambda}$ in  (\ref{Mdef}) given the restriction in \textbf{Step 1}.  
\item \textbf{Step 3}: Find an upper bound on $V(\textbf{x}_i) - \widetilde{V}(\textbf{x}_i)$ (\ref{boundO}) in terms of $h_{max}$. 
\end{description}

\noindent \textbf{Step 1}:
Define $\phi: \overline{\Omega} \rightarrow \mathbb{R}$,
\begin{equation} \label{phisubdef}\displaystyle \phi(\textbf{x}) = \frac{1}{\lambda} \left( M_{\epsilon,\lambda} + \widetilde{V}(\overline{\textbf{x}}_i) + \frac{\left\|\textbf{x} - \overline{\textbf{x}}_i\right\|^2}{2\epsilon}\right) \text{ and so } \nabla \phi(\textbf{x}) = \frac{1}{\lambda} \left(\frac{\textbf{x} - \overline{\textbf{x}}_i}{\epsilon}\right).\end{equation}

Using (\ref{phidef}), (\ref{Mdef}), (\ref{phisubdef}), and $M_{\epsilon,\lambda} \geq \Phi(\textbf{x},\overline{\textbf{x}}_i)$, it can be shown that $V(\textbf{x}) \leq \phi(\textbf{x})$ for all $\textbf{x} \in \overline{\Omega}$ and $V(\overline{\textbf{x}}) = \phi(\overline{\textbf{x}})$. Therefore $V - \phi$ has a local maximum at $\overline{\textbf{x}}$. By Lemma \ref{equivalentsub}, $\textbf{p} = \nabla \phi(\overline{\textbf{x}}) \in D^+V(\overline{\textbf{x}})$. By Lemma \ref{onesidedlemma}, $|\nabla \phi(\overline{\textbf{x}})|$ is bounded by the Lipschitz constant $L_V$, which by (\ref{overallbound}) and (\ref{phisubdef}), 
$$ \left\|\overline{\textbf{x}} - \overline{\textbf{x}}_i\right\| \leq \lambda \left\|\nabla \phi(\overline{\textbf{x}})\right\|\epsilon \leq \lambda C_0\epsilon.$$
From (\ref{lambdachoice}), and using $0 < \lambda < 1$,
\begin{equation} \label{step1}
(1-\lambda) > \frac{1}{G_{min}}\left( \frac{C_1}{\epsilon}h_{max} + \lambda L_g\left\|\overline{\textbf{x}} - \overline{\textbf{x}}_i\right\| \right).
\end{equation}
Define $\psi: \overline{\Omega}_X \rightarrow \mathbb{R}$, \begin{equation} \label{psidef} \psi(\textbf{x}_i) = -M_{\epsilon,\lambda} + \lambda V(\overline{\textbf{x}}) - \frac{\left\|\overline{\textbf{x}} - \textbf{x}_i\right\|^2}{2\epsilon}, \text{ and so } \nabla \psi(\textbf{x}_i) = \frac{\overline{\textbf{x}}-\textbf{x}_i}{\epsilon}.\end{equation}
Let $\textbf{u}_{\overline{\textbf{x}}_i}^*$ optimize the Hamiltonian (\ref{hamdef}) for arguments $\overline{\textbf{x}}_i$ and $\nabla \psi(\overline{\textbf{x}}_i)$, $$H(\overline{\textbf{x}}_i,\nabla \psi(\overline{\textbf{x}}_i)) = -\nabla \psi(\overline{\textbf{x}}_i) \cdot \textbf{u}_{\overline{\textbf{x}}_i}^* - g(\overline{\textbf{x}}_i,\textbf{u}^*_{\overline{\textbf{x}}_i}).$$ From (\ref{step1}), assumptions \textbf{(P2)}, \textbf{(P3)} and definitions of $\nabla \phi$ (\ref{phidef}) and $\nabla \psi$ (\ref{psidef}),
\begin{equation*}
(1-\lambda) g(\overline{\textbf{x}}_i,\textbf{u}_{\overline{\textbf{x}}_i}^*) > \frac{C_1}{\epsilon}{h_{max}} + \lambda( g(\overline{\textbf{x}},\textbf{u}_{\overline{\textbf{x}}_i}^*) - g(\overline{\textbf{x}}_i,\textbf{u}_{\overline{\textbf{x}}_i}^*)),
\end{equation*} 
\begin{equation*} \displaystyle \frac{\overline{\textbf{x}}-\overline{\textbf{x}}_i}{\epsilon} \cdot \textbf{u}^*_{\overline{\textbf{x}}_i} + g(\overline{\textbf{x}}_i, \textbf{u}^*_{\overline{\textbf{x}}_i})  - \lambda \left(\frac{1}{\lambda}\cdot \frac{\overline{\textbf{x}}-\overline{\textbf{x}}_i}{\epsilon} \cdot \textbf{u}^*_{\overline{\textbf{x}}_i} + g(\overline{\textbf{x}},\textbf{u}_{\overline{\textbf{x}}_i}^*)  \right) > \frac{C_1}{\epsilon}h_{max}, \end{equation*}
\begin{equation} \label{step1almostend} \displaystyle \nabla \psi(\overline{\textbf{x}}_i) \cdot \textbf{u}^*_{\overline{\textbf{x}}_i} + g(\overline{\textbf{x}}_i, \textbf{u}^*_{\overline{\textbf{x}}_i}) - \lambda(\nabla \phi(\overline{\textbf{x}}) \cdot \textbf{u}^*_{\overline{\textbf{x}}_i} + g(\overline{\textbf{x}},\textbf{u}^*_{\overline{\textbf{x}}_i})) > \frac{C_1}{\epsilon}h_{max}.\end{equation}
Since $\textbf{u}_{\overline{\textbf{x}}_i}^*$ is not necessarily the maximizer of $H(\overline{\textbf{x}},\nabla \phi(\overline{\textbf{x}}))$,
\begin{equation} \label{step1last}- \lambda(\nabla \phi(\overline{\textbf{x}}) \cdot \textbf{u}^*_{\overline{\textbf{x}}_i} + g(\overline{\textbf{x}},\textbf{u}^*_{\overline{\textbf{x}}_i})) \leq \lambda H(\overline{\textbf{x}}, \nabla \phi (\overline{\textbf{x}})). \end{equation}
It will now be shown that at most one of $\overline{\textbf{x}}_i$ or $\overline{\textbf{x}}$ can be in $\Omega$. Following (\ref{step1almostend}) and using  the definition of the Hamiltonian (\ref{hamdef}), (\ref{step1last}), $\textbf{u}_{\overline{\textbf{x}}_i}^*$ is the optimizer of $H(\overline{\textbf{x}}_i,\nabla \psi(\overline{\textbf{x}}_i))$,
\begin{equation} \label{contradictboth}  - H(\overline{\textbf{x}}_i,\nabla \psi(\overline{\textbf{x}}_i)) + \lambda H(\overline{\textbf{x}},\nabla \phi(\overline{\textbf{x}})) > \frac{C_1}{\epsilon}{h_{max}}. \end{equation}

\noindent \textbf{Case 1}: Let $\overline{\textbf{x}} \in \Omega$. From Definition \ref{subsol}, $H(\overline{\textbf{x}},\nabla \phi(\overline{\textbf{x}})) \leq 0.$ From (\ref{contradictboth}),
\begin{equation} H(\overline{\textbf{x}}_i,\nabla \psi(\overline{\textbf{x}}_i) ) < - \frac{C_1}{\epsilon}h_{max}. \label{useproof} \end{equation}
For all $\textbf{x}_i \in X_0$, $\psi(\textbf{x}_i) \leq \widetilde{V}(\textbf{x}_i)$, $\psi(\overline{\textbf{x}}_i) = \widetilde{V}(\overline{\textbf{x}}_i)$. By Definition \ref{dcstencil} and Theorem \ref{monotonic},
\begin{equation} \label{monotonicoumproof} \widetilde{H}[\overline{\textbf{NF}}(\textbf{x}_i),\widetilde{V}](\overline{\textbf{x}}_i,\widetilde{V}(\overline{\textbf{x}}_i)) = \widetilde{H}[S(\textbf{x}_i),\widetilde{V}](\overline{\textbf{x}}_i,\widetilde{V}(\overline{\textbf{x}}_i)) \leq \widetilde{H}[S(\textbf{x}_i),\psi](\overline{\textbf{x}}_i,\psi(\overline{\textbf{x}}_i)).\end{equation}
It will be shown that $\overline{\textbf{x}}_i \in X_0\cap\Omega^c$ using proof by contrapositive. Since $\widetilde{V}$ is a solution to the numerical HJB equation (\ref{numhamstatichjb}), it is a supersolution of the numerical HJB equation (Definition \ref{supnumham}). If $\overline{\textbf{x}}_i \in X_0\cap\Omega$, then
\begin{equation} \widetilde{H}[\overline{\textbf{NF}}(\textbf{x}_i),\widetilde{V}](\overline{\textbf{x}}_i,\widetilde{V}(\overline{\textbf{x}}_i)) = \widetilde{H}[S(\textbf{x}_i),\widetilde{V}](\overline{\textbf{x}}_i,\widetilde{V}(\overline{\textbf{x}}_i))  \geq 0. \label{prop1} \end{equation}
Furthermore if $\overline{\textbf{x}}_i \in X_0 \cap \Omega$, Theorem \ref{consistent} must also hold. That is, since $\left\|\nabla^2 \psi \right\|_2 = \frac{1}{\epsilon}$,
\begin{equation} \label{prop2} |H(\overline{\textbf{x}}_i,\nabla \psi(\overline{\textbf{x}}_i)) - \widetilde{H}[S(\textbf{x}_i),\psi](\overline{\textbf{x}}_i,\psi(\overline{\textbf{x}}_i))| \leq \frac{C_1}{\epsilon}h_{max}. \end{equation}
It will be shown (\ref{prop1}) and (\ref{prop2}) cannot simultaneously be true, implying $\overline{\textbf{x}}_i \in X_0 \cap \Omega^c$. If (\ref{prop1}) is true, then by (\ref{monotonicoumproof}), $\widetilde{H}[S(\textbf{x}_i),\psi](\overline{\textbf{x}}_i,\psi(\overline{\textbf{x}}_i)) \geq 0$. By (\ref{useproof}),
$$ H(\overline{\textbf{x}}_i, \nabla \psi(\overline{\textbf{x}}_i)) - \widetilde{H}[S(\textbf{x}_i),\psi](\overline{\textbf{x}}_i,\psi(\overline{\textbf{x}}_i)) < -\frac{C_1}{\epsilon}h_{max}.$$
Therefore (\ref{prop2}) is false. 

Otherwise, if (\ref{prop2}) were true, using (\ref{useproof}), $$ H(\overline{\textbf{x}}_i, \nabla \psi(\overline{\textbf{x}}_i)) - \widetilde{H}[S(\textbf{x}_i),\psi](\overline{\textbf{x}}_i,\psi(\overline{\textbf{x}}_i)) \geq -\frac{C_1}{\epsilon} h_{max} > H(\overline{\textbf{x}}_i, \nabla \psi(\overline{\textbf{x}}_i)). $$
Hence with (\ref{monotonicoumproof}),
$$\widetilde{H}[\overline{\textbf{NF}}(\textbf{x}_i),\widetilde{V}](\overline{\textbf{x}}_i,\widetilde{V}(\overline{\textbf{x}}_i)) =\widetilde{H}[S(\textbf{x}_i),\widetilde{V}](\overline{\textbf{x}}_i,\widetilde{V}(\overline{\textbf{x}}_i)) \leq \widetilde{H}[S(\textbf{x}_i),\psi](\overline{\textbf{x}}_i,\psi(\overline{\textbf{x}}_i)) < 0.$$ Therefore (\ref{prop1}) is false. Hence $\overline{\textbf{x}}_i \in X_0 \cap \Omega^c$.

\noindent \textbf{Case 2}: If $\overline{\textbf{x}}_i \in X_0 \cap \Omega$, from Theorem \ref{consistent},
\begin{equation} \label{step12} \widetilde{H}[S(\textbf{x}_i),\psi](\overline{\textbf{x}}_i,\psi(\textbf{x}_i)) - H(\overline{\textbf{x}}_i,\nabla \psi(\overline{\textbf{x}}_i)) \leq  \frac{C_1}{\epsilon}h_{max}.\end{equation}
From (\ref{monotonicoumproof}), Definition \ref{dcstencil} and $\widetilde{V}$ is a supersolution of the numerical HJB (\ref{numhamstatichjb}) (Definition \ref{supnumham}),
 \begin{equation*} \widetilde{H}[S(\textbf{x}_i),\psi](\overline{\textbf{x}}_i, \psi(\overline{\textbf{x}}_i)) \geq \widetilde{H}[S(\textbf{x}_i),\widetilde{V}](\overline{\textbf{x}}_i, \widetilde{V}(\overline{\textbf{x}}_i)) = \widetilde{H}[\overline{\textbf{NF}}(\textbf{x}_i),\widetilde{V}](\overline{\textbf{x}}_i, \widetilde{V}(\overline{\textbf{x}}_i)) \geq 0. \end{equation*}
From (\ref{contradictboth}) and (\ref{step12}),
\begin{equation} \label{step122} \frac{C_1}{\epsilon}h_{max} + \widetilde{H}[S(\textbf{x}_i),\psi](\overline{\textbf{x}}_i,\psi(\overline{\textbf{x}}_i)) -\lambda H(\overline{\textbf{x}},\nabla \phi(\overline{\textbf{x}})) < \frac{C_1}{\epsilon} h_{max}. \end{equation}
Since $\overline{\textbf{x}}_i \in X_0 \cap \Omega$, $\widetilde{H}[S(\textbf{x}_i),\psi](\overline{\textbf{x}}_i,\psi(\overline{\textbf{x}}_i)) \geq 0$, from (\ref{step122}), and $0<\lambda<1$, $$ H(\overline{\textbf{x}},\nabla \phi(\overline{\textbf{x}})) > 0,$$ which implies by Definition \ref{subsol} of the viscosity subsolution, $\overline{\textbf{x}} \in \partial \Omega$.  Hence at most one of maximizers of $M_{\epsilon,\lambda}$, $\overline{\textbf{x}}$ and $\overline{\textbf{x}}_i$ can belong to $\Omega$.

\noindent \textbf{Step 2}: An upper bound on $M_{\epsilon,\lambda}$ (\ref{Mdef}) will be found. 	

\noindent \textbf{Case 1}:  $\overline{\textbf{x}} \in \overline{\Omega}$, $\overline{\textbf{x}}_i \in X_0 \cap \Omega^c$.

Let $\check{\textbf{x}} = \argmin_{\textbf{x} \in \partial \Omega} \left\|\overline{\textbf{x}}_i - \textbf{x}\right\|$. Let $\underline{\textbf{x}}_i$ be the point on the line from $\overline{\textbf{x}}$ and $\overline{\textbf{x}}_i$ intersecting $\partial\Omega$ . For $\overline{\textbf{x}} \in \partial \Omega$, $\underline{\textbf{x}}_i = \overline{\textbf{x}}$. Since $\overline{\Omega}$ is convex, by Lemma \ref{obtuseconvex},  the angle between vectors $\underline{\textbf{x}}_i - \check{\textbf{x}}$ and $\overline{\textbf{x}}_i - \check{\textbf{x}}$ is nonacute. Using the cosine law,
\begin{align*}
\left\|\underline{\textbf{x}}_i - \overline{\textbf{x}}_i\right\|^2 &= \left\|\underline{\textbf{x}}_i - \check{\textbf{x}}\right\|^2 + \left\|\overline{\textbf{x}}_i - \check{\textbf{x}}\right\|^2 - 2  ( \underline{\textbf{x}}_i - \check{\textbf{x}}) \cdot( \overline{\textbf{x}}_i -\check{\textbf{x}}), \\
&\geq \left\|\underline{\textbf{x}}_i - \check{\textbf{x}}\right\|^2,\\
\left\|\underline{\textbf{x}}_i - \overline{\textbf{x}}_i\right\| & \geq  \left\|\underline{\textbf{x}}_i - \check{\textbf{x}}\right\|.
\end{align*}
Since $\underline{\textbf{x}}_i$ is on the line segment from $\overline{ \textbf{x}}$ to $\overline{\textbf{x}}_i$, $\left\|\overline{\textbf{x}} - \overline{\textbf{x}}_i\right\|= \left\|\overline{\textbf{x}} - \underline{\textbf{x}}_i\right\| + \left\|\underline{\textbf{x}}_i - \overline{\textbf{x}}_i\right\|$. With the triangle inequality,
\begin{align}
\left\|\overline{\textbf{x}} - \underline{\textbf{x}}_i\right\| + \left\|\underline{\textbf{x}}_i - \overline{\textbf{x}}_i\right\| & \geq \left\|\overline{\textbf{x}} - \underline{\textbf{x}}_i\right\| + \left\|\underline{\textbf{x}}_i - \check{\textbf{x}}\right\|, \nonumber \\
\left\|\overline{\textbf{x}} - \overline{\textbf{x}}_i\right\| & \geq \left\|\overline{\textbf{x}} - \check{\textbf{x}}\right\|. \label{triangleineq}
\end{align}
By the Lipschitz-continuity of $V$ with constant $C_0$, $0 < \lambda < 1$, $|\widetilde{V}| \leq C_0$, and since $\widetilde{V}$ is a supersolution to the numerical HJB equation (\ref{numhamstatichjb}), $\widetilde{V}(\overline{\textbf{x}}_i) \geq q(\check{\textbf{x}})$, 
\begin{align}
\nonumber M_{\epsilon,\lambda} & = \lambda V(\overline{\textbf{x}}) - \widetilde{V}(\overline{\textbf{x}}_i) - \frac{\left\|\overline{\textbf{x}} - \overline{\textbf{x}}_i\right\|^2}{2\epsilon},\\
\nonumber & = \lambda( V(\overline{\textbf{x}}) - \widetilde{V}(\overline{\textbf{x}}_i)) - (1-\lambda)\widetilde{V}(\overline{\textbf{x}}_i) - \frac{\left\|\overline{\textbf{x}} - \overline{\textbf{x}}_i\right\|^2}{2\epsilon},\\
& \leq  \lambda( V(\overline{\textbf{x}}) - q(\check{\textbf{x}})) + (1-\lambda)C_0 - \frac{\left\|\overline{\textbf{x}} - \overline{\textbf{x}}_i\right\|^2}{2\epsilon} \label{label},
\end{align}
If $\overline{\textbf{x}} \in \Omega$, $ V(\check{\textbf{x}}) \leq q(\check{\textbf{x}})$, from (\ref{label}),
\begin{equation} \label{xinomega}
M_{\epsilon,\lambda} \leq \lambda( V(\overline{\textbf{x}}) -  V(\check{\textbf{x}})) +(1-\lambda)C_0 - \frac{\left\|\overline{\textbf{x}} - \overline{\textbf{x}}_i\right\|^2}{2\epsilon}.
\end{equation}
Otherwise $\overline{\textbf{x}} \in \partial \Omega$, and $V(\overline{\textbf{x}}) \leq q(\overline{\textbf{x}})$, from (\ref{label}),
\begin{equation} \label{xonboundary}
M_{\epsilon,\lambda}  \leq  \lambda(q(\overline{\textbf{x}}) - q(\check{\textbf{x}})) + (1-\lambda)C_0  - \frac{\left\|\overline{\textbf{x}} - \overline{\textbf{x}}_i\right\|^2}{2\epsilon}.\\
\end{equation}
The Lipschitz continuity of both $q$ and $V$ with constant $C_0$ in (\ref{xinomega}) and (\ref{xonboundary}) and $\left\|\overline{\textbf{x}} - \overline{\textbf{x}}_i\right\| \geq \left\|\overline{\textbf{x}} - \check{\textbf{x}}\right\|$ from (\ref{triangleineq}), along with $0<\lambda<1$ yield
\begin{equation}
M_{\epsilon,\lambda} \leq  C_0\left\|\overline{\textbf{x}} - \check{\textbf{x}}\right\| + (1-\lambda)C_0  - \frac{\left\|\overline{\textbf{x}} - \check{\textbf{x}}\right\|^2}{2\epsilon},\\
\end{equation}
which is quadratic in $\left\|\overline{\textbf{x}} - \check{\textbf{x}}\right\|$. The quadratic is maximized with $\left\|\overline{\textbf{x}} - \check{\textbf{x}}\right\| = C_0\epsilon$. Thus,
\begin{equation} \label{case2bound} M_{\epsilon,\lambda} \leq (1-\lambda)C_0 + \frac{C_0^2 \epsilon}{2}.\end{equation}

\noindent \textbf{Case 2}: $\overline{\textbf{x}} \in \partial \Omega$, $\overline{\textbf{x}}_i \in X_0 \cap \Omega$.

From Lemma \ref{boundarylemma}, there exists $\hat{\textbf{x}}_i \in  X_0 \cap \Omega^c$ such that
\begin{equation} \label{closepoint}  \left\|\overline{\textbf{x}} - \hat{\textbf{x}}_i\right\| \leq h_{max}. \end{equation}
Let $\widetilde{\textbf{x}} = \argmin_{\textbf{x} \in \partial \Omega} \left\|\hat{\textbf{x}}_i - \textbf{x}\right\|$. Using $0 < \lambda < 1$, $\widetilde{V}(\hat{\textbf{x}}_i) \geq q(\widetilde{\textbf{x}})$, 	$ V(\overline{\textbf{x}}) \leq q(\overline{\textbf{x}})$, Lipschitz-continuity of $q$ and $\widetilde{V}$ both with constant $C_0$,
\begin{align*}
M_{\epsilon,\lambda} & = \lambda V(\overline{\textbf{x}}) - \widetilde{V}(\overline{\textbf{x}}_i) - \frac{\left\|\overline{\textbf{x}} - \overline{\textbf{x}}_i\right\|^2}{2\epsilon},\\
& =  \lambda( V(\overline{\textbf{x}}) - \widetilde{V}(\overline{\textbf{x}}_i)) - (1-\lambda)\widetilde{V}(\overline{\textbf{x}}_i) - \frac{\left\|\overline{\textbf{x}} - \overline{\textbf{x}}_i\right\|^2}{2\epsilon},\\
& \leq \lambda(q(\overline{\textbf{x}}) - q(\widetilde{\textbf{x}}) + q(\widetilde{\textbf{x}}) - \widetilde{V}(\overline{\textbf{x}}_i)) + (1-\lambda)C_0 - \frac{\left\|\overline{\textbf{x}} - \overline{\textbf{x}}_i\right\|^2}{2\epsilon}, \\
& \leq \lambda C_0\left\|\overline{\textbf{x}} - \widetilde{\textbf{x}}\right\| + \lambda(\widetilde{V}(\hat{\textbf{x}}_i) - \widetilde{V}(\overline{\textbf{x}}_i)) + (1-\lambda)C_0 - \frac{\left\|\overline{\textbf{x}} - \overline{\textbf{x}}_i\right\|^2}{2\epsilon}, \\
& \leq  C_0 (\left\|\overline{\textbf{x}} - \widetilde{\textbf{x}}\right\| + \left\|\hat{\textbf{x}}_i - \overline{\textbf{x}}_i\right\|) + (1-\lambda)C_0 - \frac{\left\|\overline{\textbf{x}} - \overline{\textbf{x}}_i\right\|^2}{2\epsilon},
\end{align*}
Using the triangle inequality, $\left\|\hat{\textbf{x}}_i - \overline{\textbf{x}}_i\right\| \leq  \left\|\hat{\textbf{x}}_i - \widetilde{\textbf{x}}\right\| + \left\|\widetilde{\textbf{x}} - \overline{\textbf{x}}\right\| + \left\|\overline{\textbf{x}} - \overline{\textbf{x}}_i\right\|$, hence
\begin{equation*}
M_{\epsilon,\lambda} \leq (1-\lambda)C_0 +  C_0 (\left\|\overline{\textbf{x}} - \widetilde{\textbf{x}}\right\| + \left\|\hat{\textbf{x}}_i - \widetilde{\textbf{x}}\right\| + \left\|\widetilde{\textbf{x}} - \overline{\textbf{x}}\right\| + \left\|\overline{\textbf{x}} - \overline{\textbf{x}}_i\right\|) - \frac{\left\|\overline{\textbf{x}} - \overline{\textbf{x}}_i\right\|^2}{2\epsilon}. \\ 
\end{equation*}
By Lemma \ref{obtuseconvex}, and the cosine law, $\left\|\overline{\textbf{x}} - \widetilde{\textbf{x}}\right\| \leq \left\|\overline{\textbf{x}} - \hat{\textbf{x}}_i\right\|$. From the definition of $\widetilde{\textbf{x}}$, $\left\|\hat{\textbf{x}}_i-\widetilde{\textbf{x}}\right\| \leq \left\|\overline{\textbf{x}} - \hat{\textbf{x}}_i\right\|$. Therefore, $$ M_{\epsilon,\lambda} \leq (1-\lambda)C_0 + 3C_0\left\|\overline{\textbf{x}} - \hat{\textbf{x}}_i\right\| + C_0\left\|\overline{\textbf{x}} - \overline{\textbf{x}}_i\right\| - \frac{\left\|\overline{\textbf{x}} - \overline{\textbf{x}}_i\right\|^2}{2\epsilon}.$$
From (\ref{closepoint}) and maximizing over the quadratic $\left\|\overline{\textbf{x}} - \overline{\textbf{x}}_i\right\|$ with $\left\|\overline{\textbf{x}} - \overline{\textbf{x}}_i\right\| = C_0\epsilon$, 
\begin{equation}
M_{\epsilon,\lambda} \leq (1-\lambda)C_0 + 3C_0h_{max} + \frac{C_0^2\epsilon}{2}. \label{case3bound}
\end{equation}
\noindent \textbf{Step 3}: The upper bound of $M_{\epsilon,\lambda}$ in (\ref{case3bound}) is larger than (\ref{case2bound}). From (\ref{boundO}),
\begin{align*} \displaystyle
V(\textbf{x}_i) - \widetilde{V}(\textbf{x}_i) & \leq  C_0(1 - \lambda) + M_{\epsilon,\lambda}, \\
& \leq 2C_0(1-\lambda) + 3C_0h_{max} + \frac{C_0^2 \epsilon}{2},\\
& \leq 2C_0\frac{2}{G_{min}}\left(\frac{C_1}{\epsilon}h_{max} + C_0L_g\epsilon \right) + 3C_0h_{max} + \frac{C_0^2 \epsilon}{2}, \\
& \leq \left( \frac{4C_0C_1}{G_{min}} + \frac{4C_0^2 L_g}{G_{min}} + \frac{C_0^2}{2} \right) \left(\frac{h_{max}}{\epsilon} + \epsilon\right) + 3C_0h_{max},
\end{align*}
Since $\epsilon = \sqrt{h_{max}}$  is a global minimum of $(\frac{h_{max}}{\epsilon} + \epsilon)$, and setting \\
$C = 2\left( \frac{4C_0C_1}{G_{min}} + \frac{4C_0^2 L_g}{G_{min}} + \frac{C_0^2}{2} + 3C_0 \right)$, for $h_{max} <h_0 = \min\{\frac{G_{min}^2}{4(C_1+C_0L_g)^2},1\}$,
\begin{equation}
V(\textbf{x}_i) - \widetilde{V}(\textbf{x}_i) \leq  C \sqrt{h_{max}}.
\label{finalresult}
\end{equation}
Finally, a symmetrical argument using $V$ a viscosity supersolution of (\ref{conttheo1}) (Definition \ref{supersol}), and $\widetilde{V}$ a subsolution of the numerical HJB equation (\ref{numhamstatichjb}) (Definition \ref{subnumham}) can show (\ref{finalresult}) with $V(\textbf{x}_i)$ and $\widetilde{V}(\textbf{x}_i)$ reversed. Hence for $h_{max} <  h_0$,
$$ \max_{\textbf{x}_i \in X_0} |\widetilde{V}(\textbf{x}_i) - V(\textbf{x}_i)| \leq C\sqrt{h_{max}}. \square$$
Theorem \ref{oumerror} will now be extended to $\overline{\Omega}_X$. Define $\hat{V}:\overline{\Omega}_X \rightarrow \mathbb{R}$
$$\hat{V}(\textbf{x}) = \sum_{j=0}^n \zeta_j V(\textbf{x}_j^{\textbf{s}}) \text{ for } \textbf{x} = \sum_{j=0}^n \zeta_j \textbf{x}_j^{\textbf{s}}.$$ On $\textbf{x}_i \in X_0$, $V(\textbf{x}_i) = \hat{V}(\textbf{x}_i)$ are equal.
\begin{lemma} \label{otherfunction}
There exists $D_1>0$ for all $\textbf{x} \in \overline{\Omega}_X$, such that
\begin{equation} |V(\textbf{x}) - \hat{V}(\textbf{x})| \leq D_1h_{max}.\end{equation}\end{lemma}
\noindent \textit{Proof}. Let $\zeta \in \Xi_n$ and $\textbf{x} \in \textbf{s}$ such that $\textbf{x} = \sum_{j=0}^n \zeta_j \textbf{x}_j^{\textbf{s}}$. Using $V(\textbf{x}_i) = \hat{V}(\textbf{x}_i)$ for all vertices $\textbf{x}_i \in X_0$, $\sum_{j=0}^n \zeta_j = 1$, Lemma \ref{valuelipschitz}, with Lipschitz constant $L_V=2G_{max}$,
\begin{equation*}\displaystyle
|V(\textbf{x}) - \hat{V}(\textbf{x})| \leq \sum_{j=0}^n \zeta_j |V(\textbf{x}) - V(\textbf{x}_j)| \leq 2G_{max}h_{max}. \square
\end{equation*}
\begin{corollary}
There exists $D_2>0$ for all $\textbf{x} \in \overline{\Omega}_X$ such that 
\begin{equation} |V(\textbf{x}) - \widetilde{V}(\textbf{x})| \leq D_2\sqrt{h_{max}},\end{equation}
for $h_{max} < h_0$ as described in Theorem \ref{oumerror}.
\end{corollary}

\noindent \textit{Proof.}
Let $\zeta \in \Xi_n$ and $\textbf{x} \in \textbf{s}$ such that $\textbf{x} = \sum_{j=0}^n \zeta_j \textbf{x}_j^{\textbf{s}}$. For $\textbf{x} \in \overline{\Omega}_X$, $\widetilde{V}(\textbf{x}) = \sum_{j=0}^n \zeta_j \widetilde{V}(\textbf{x}_j^{\textbf{s}})$. From Lemma \ref{otherfunction} and Theorem \ref{oumerror},
\begin{align*}
V(\textbf{x}) - \widetilde{V}(\textbf{x}) & \leq D_1h_{max} + \hat{V}(\textbf{x}) - \widetilde{V}(\textbf{x})\\
& = D_1h_{max} +  \sum_{j=0}^n \zeta_j (\hat{V}(\textbf{x}_j^{\textbf{s}}) - \widetilde{V}(\textbf{x}_j^{\textbf{s}}))\\
& \leq (D_1+C)\sqrt{h_{max}},
\end{align*}
for $h_{max} < h_0$. The proof for $\widetilde{V}(\textbf{x}) - V(\textbf{x})$ is symmetrical. Hence $D_2 = D_1 + C$.$\square$

\section{Numerical Convergence of OUM Example} \label{oumconverge}

\begin{figure}
\begin{center} \subfloat[][Rectangular speed profile $\mathcal{U}_g(\textbf{x})$ with length $6$ in the $x$-direction and $2$ in the $y$-direction.]{\includegraphics[width=2in]{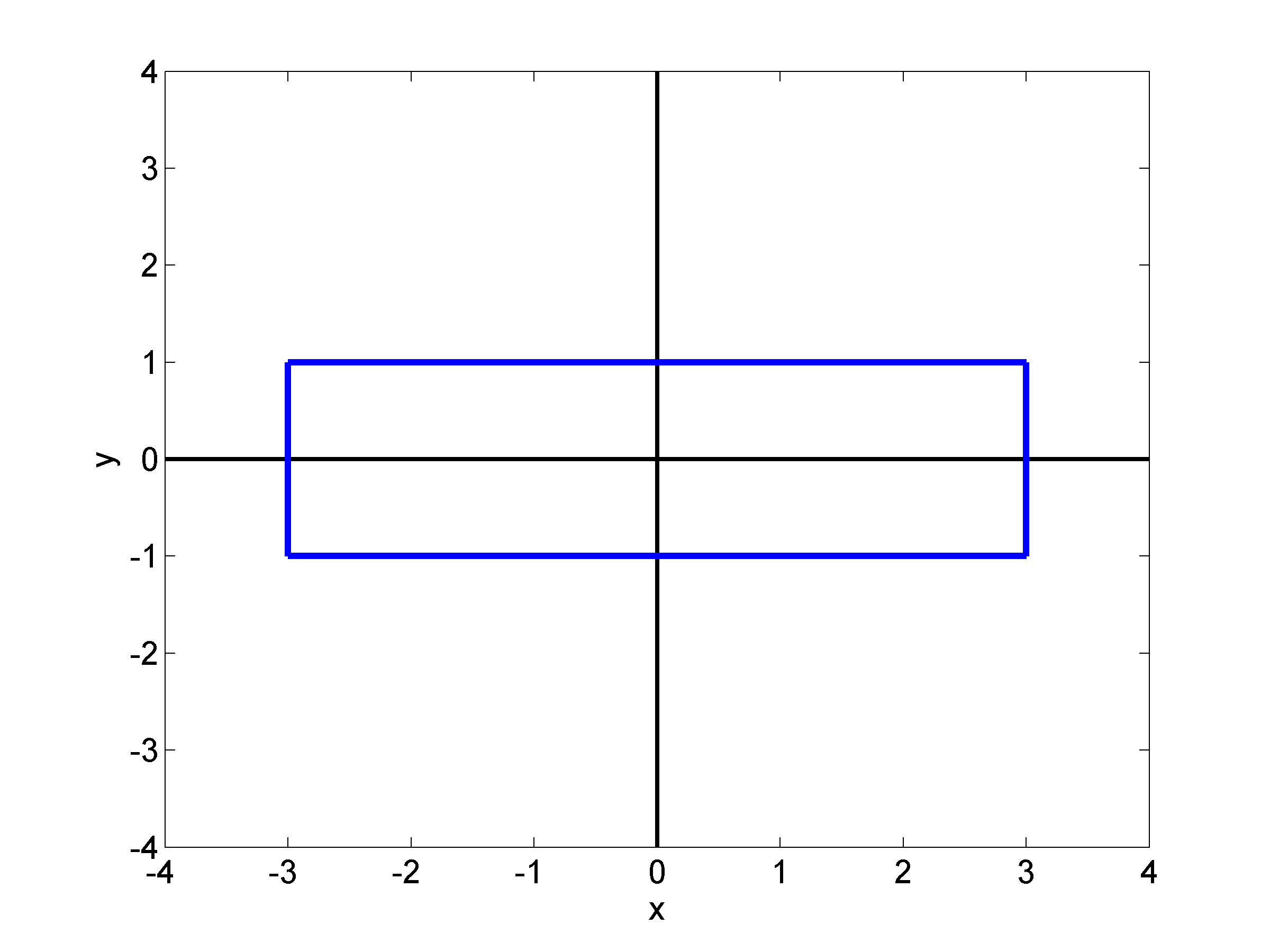} \label{rectspeedprofile}}\hspace{2mm}
\subfloat[][The exact solution $V$: a three-dimensional view.]{\includegraphics[width=2in]{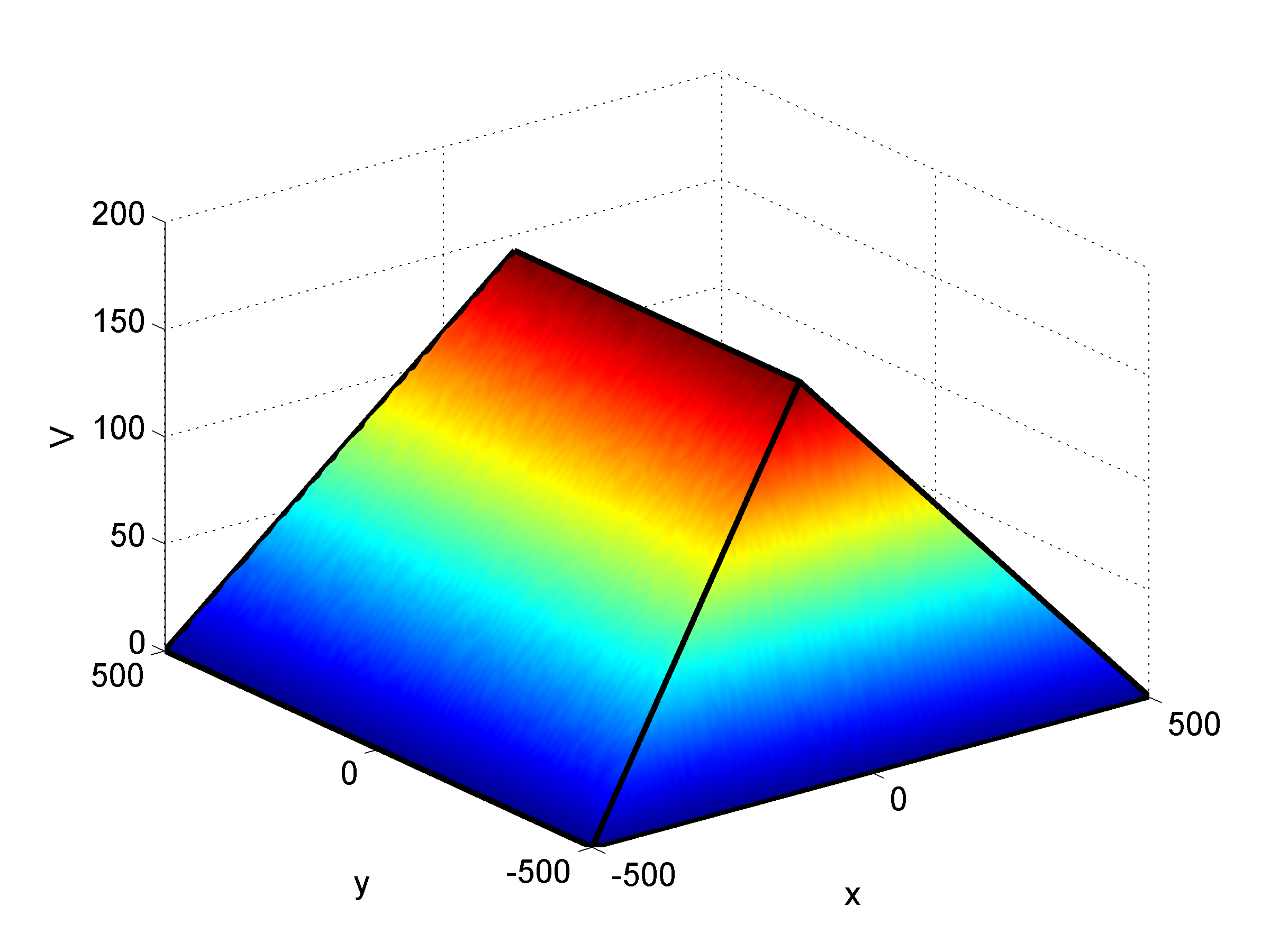} \label{exactsolution3}} \hspace{2mm} \subfloat[][The error between $\widetilde{V}$ and $V$ (viewed from above) is greatest at points where $\nabla V$ is not defined.]{\includegraphics[width=2in]{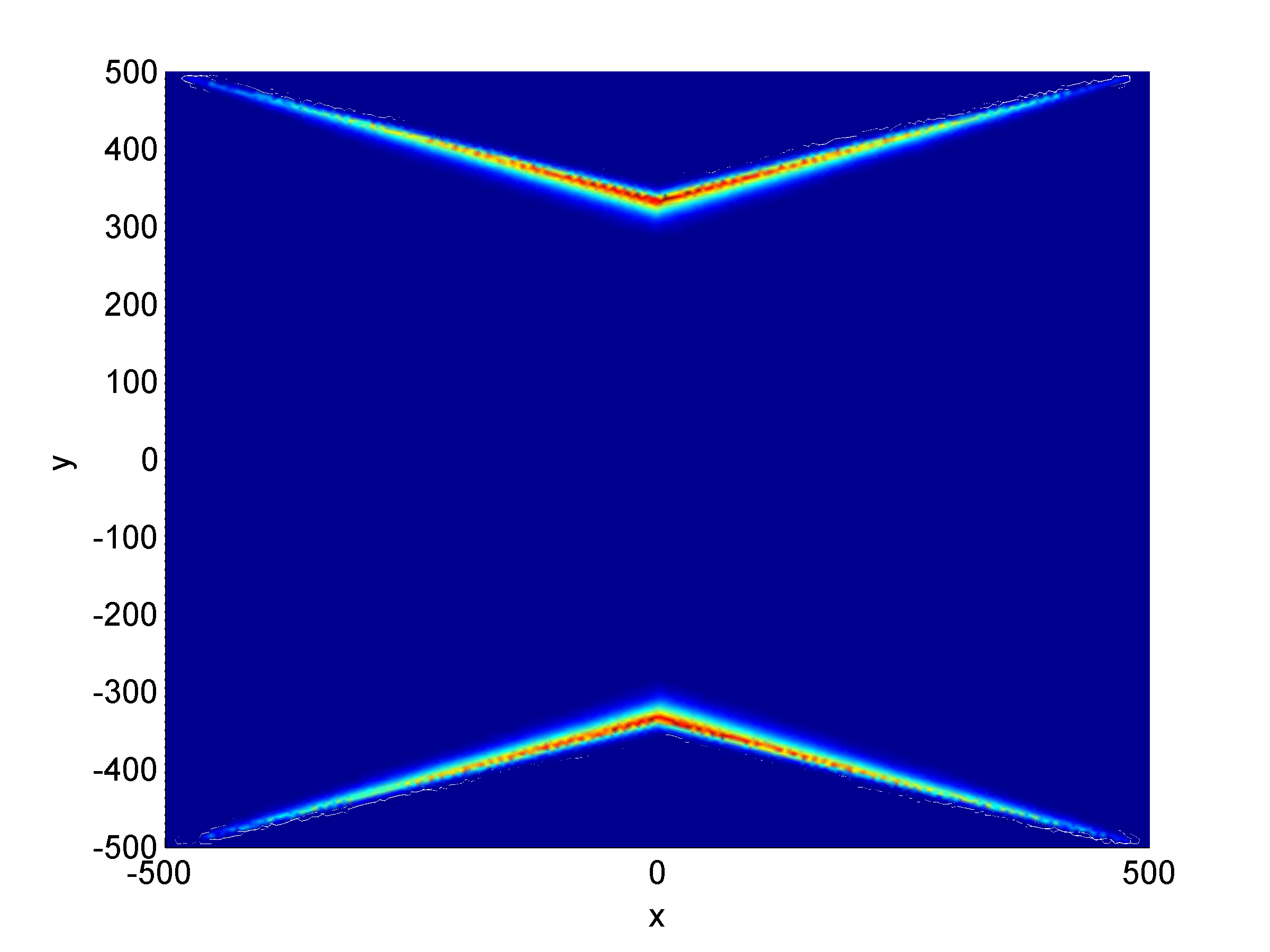} \label{errorrectplot}} 
 \end{center}
\caption{A numerical example: a) speed profile, b) true solution $V$, c) error.}
\label{exactsolution}
\end{figure}

An example of the error computed using OUM for the boundary value problem is given. The OUM algorithm was programmed in MATLAB$^{\text{\textregistered}}$ on an ASUS X550L Laptop with Intel$^\text{\textregistered}$ Core $^\text{TM}$ i5 -4210U CPU Processor (1.7 GHz/2.4GHz) with 4GB RAM. As in \cite{oum}, the update for the OUM algorithm (\ref{oumupdate}) was solved using the golden section search. For $\overline{\Omega} = [-500,500]\times [-500,500]$, $\partial \Omega = \{(x,y)\in \overline{\Omega}| |x| = 500 \text{ or } |y| = 500\}$, the weight $g$ used corresponded to a rectangular speed profile (Definition \ref{speedprofiledef}) centred about \textbf{x} with dimensions $6$ in the $x$-direction and $2$ in the $y$-direction. See Figure \ref{rectspeedprofile}. The boundary function was $q(\textbf{x}) = 0$ for $\textbf{x} \in \partial \Omega$.  The same speed profile was used for all $\textbf{x} \in \Omega$. The analytic solution is made up of the concatenation of 4 planes: $y + z = 500$, $x + 3z = 500$, $-y + z = 500$ and $-x + 3z = 500$ within $\Omega$. See Figure \ref{exactsolution3}. 

\begin{figure}\begin{center}
\includegraphics[width=2.5in]{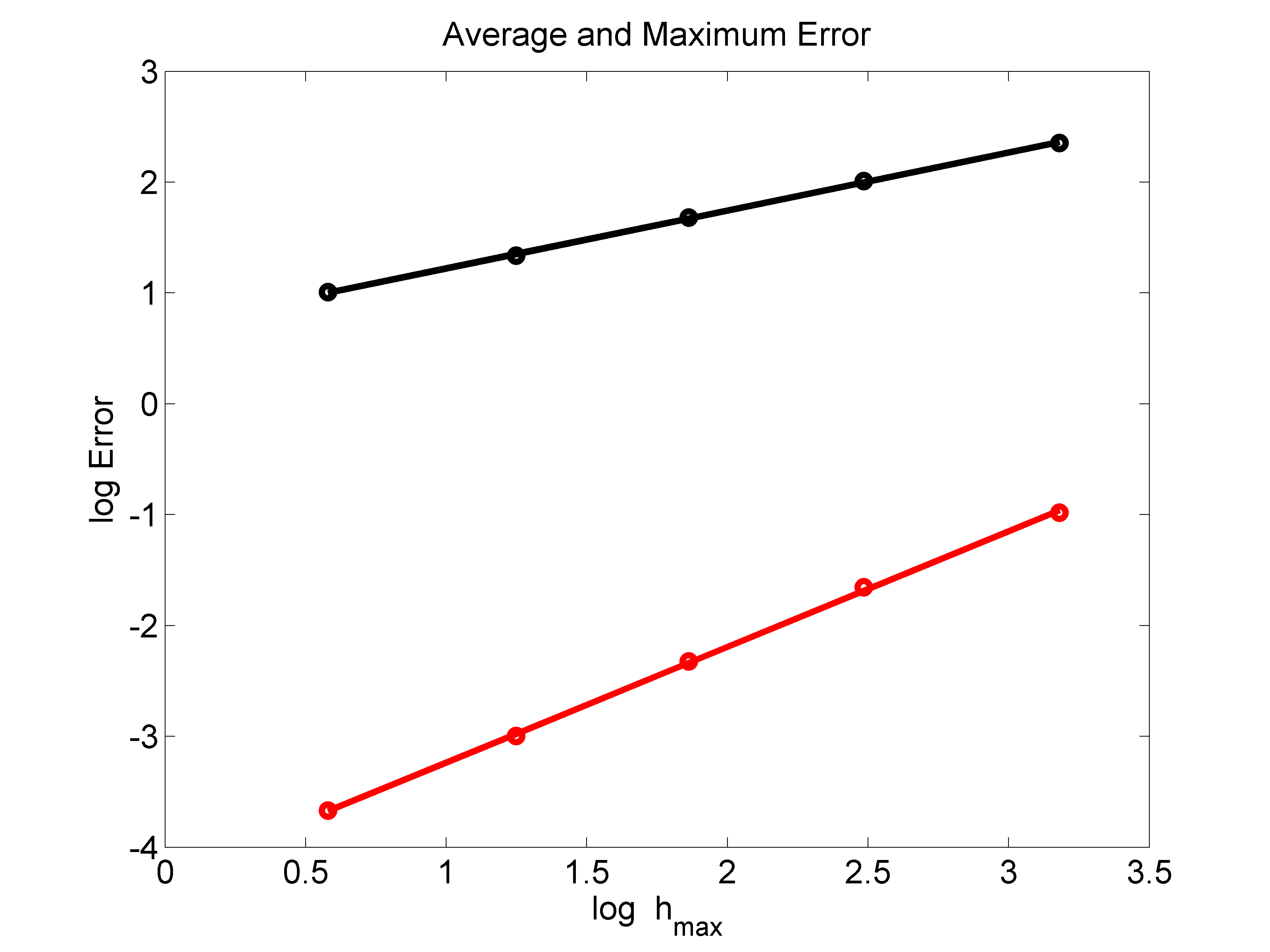} \end{center}
\caption{Average and maximum error for OUM Convergence Example - average error shown in red (below), maximum error shown in black (above). The overall convergence rates measured were $r_{avg} = 1.043$ and $r_{max} = 0.523$.}
\label{convergenceplot}
\end{figure}

\begin{table} \begin{center}
\begin{tabular}{|c|c|c|c|c|c|c|}\hline
 Vertices & Triangles& $h_{max}$ & Avg Error & $r_{avg}$ & Max Error & $r_{max}$ \\ \hline
4289 & 8256 & 24.07 & 0.3746 & -  & 10.54 & - \\ \hline
16765 & 32888 & 11.99 & 0.1914 & 0.9634 & 7.48 & 0.4931 \\ \hline
66291 & 131300 & 6.438 & 0.0979 & 1.0779 & 5.38 & 0.5289 \\ \hline
263597& 524632 & 3.483 & 0.0499 & 1.0968 & 3.80 & 0.5643 \\ \hline
1051261 & 2097400 & 1.785& 0.0255 & 1.0062 & 2.74 & 0.4900 \\ \hline
\end{tabular} \end{center}
\caption{Accuracy of OUM for a Boundary Value Problem - The OUM was used to solve the static HJB problem with a rectangular profile on five meshes. Both average error across the vertices and maximum vertex error are reported. The incremental rates of convergence are also shown.}
\label{oumaccuracy}
\end{table}
Given a set of boundary points, meshes with uneven triangles were generated using Mesh2D \cite{mesh2d}. The error values are given in Table \ref{oumaccuracy} and a plot is provided in Figure \ref{convergenceplot}. Using \texttt{polyfit} in MATLAB with the data provided in Table \ref{oumaccuracy}, affine approximations of the log-log slope fit using least squares were found. Using all 5 data points, overall rates of convergence of $r_{avg} = 1.043$ and $r_{max} = 0.523$ were obtained for average error and maximum error across the vertices respectively. The convergence rate for maximum error in this example matches closely to the theoretical results shown earlier. In average error, the OUM algorithm is at most first-order accurate (as described in \cite{oum}) since the update formula (\ref{oumupdate}) is a first-order approximation. Since $V$ is Lipschitz continuous, from Rademacher's theorem, $\nabla V$ can only be undefined on a set of measure zero. The error for all discretiztaions had the same general shape, appearing greatest near where $\nabla  V$ was undefined. See Figure \ref{errorrectplot}. Characteristics flow into, but not out of such points where $\nabla V$ is undefined, preventing the error from being propagated further \cite{kumarvlad}, hence the expected first-order convergence rate in average error.

\section{Conclusions and Future Work} \label{secconcandfw}

It was proven in this paper that the rate of convergence of the approximate solution provided by OUM to the viscosity solution of the HJB for prescribed boundary values is at least $\mathcal{O}(\sqrt{h_{max}})$ in maximum error. The basic idea of the proof is an extension of a similar proof for FMM in \cite{monneau}. A key step was to show the existence of a directionally complete stencil. This implied from existing results that the numerical Hamiltonian for the OUM is both consistent and monotonic. 

An extension of this work would be to provide a convergence rate proof for OUM in the single-source point formulation of the static HJB. This will extend the applicability of the result shown here to point-to-point path planning problems, such as for rovers \cite{shumpaper} and other robots \cite{thepathtoefficiency}. Constructing a directionally complete stencil as done here may be difficult near the source point.  

Another direction of research could be to prove that the convergence in average error of OUM is at a rate of $\mathcal{O}(h_{max})$ as was the case in the example in this paper. This could follow because OUM is a first-order method, with $V$ generally not differentiable only on a set of measure zero. Additional assumptions of regularity, such as a continuously differentiable speed profile, may lead to a proof for first-order convergence in average error applicable to many problems. 

\bibliographystyle{spmpsci}      

\begin{thebibliography}{10}
\providecommand{\url}[1]{{#1}}
\providecommand{\urlprefix}{URL }
\expandafter\ifx\csname urlstyle\endcsname\relax
  \providecommand{\doi}[1]{DOI~\discretionary{}{}{}#1}\else
  \providecommand{\doi}{DOI~\discretionary{}{}{}\begingroup
  \urlstyle{rm}\Url}\fi

\bibitem{kaltonthesis}
Alton, K.: {D}ijkstra-like ordered upwind methods for solving static
  {H}amilton-{J}acobi equations.
\newblock Ph.D. thesis, University of British Columbia (2010)

\bibitem{kaltonpaper}
Alton, K., Mitchell, I.: An ordered upwind method with precomputed stencil and
  monotone node acceptance for solving static convex {H}amilton-{J}acobi
  equations.
\newblock Journal of Scientific Computing \textbf{51}, 313--348 (2012)

\bibitem{angellnotes}
Angell, T.: Notes on convex sets.
\newblock http://www.math.udel.edu/\textasciitilde angell/Opt/convex.pdf (2011). Accessed 5 August 2013.

\bibitem{abgrall}
Augoula, S., Abgrall R.: High order numerical discretization for {H}amilton-{J}acobi equations on triangular meshes.
\newblock Journal of Scientific Computing. \textbf{15}, 197--229 (2000)

\bibitem{bardi}
Bardi, M., Capuzzo-Dolcetta, I.: Optimal Control and Viscosity Solutions of
  {H}amilton-{J}acobi-{B}ellman Equations.
\newblock Birkh\"{a}user (1997)

\bibitem{bardifalconeerror}
Bardi, M., Falcone, M.: Discrete approximation of the minimal time function for
  systems with regular optimal trajectories.
\newblock Analysis and Optimation of Systems: Lecture Notes in Control and
  Information Sciences \textbf{144}, 103--112 (1990)

\bibitem{borweinlewis}
Borwein, J.M., Lewis, A.S.: Convex Analysis and Nonlinear Optimization: Theory
  and Examples.
\newblock Springer (2005)

\bibitem{oumreactivefluxes}
Cameron, M.K.: Estimation of reactive fluxes in gradient stochastic systems
  using an analogy with electric circuits.
\newblock Journal of Computational Physics \textbf{247}, 137--152 (2013)

\bibitem{crandalllionstwoapprox}
Crandall, M.G., Lions, P.L.: Two approximations of solutions of {H}amilton-{J}acobi equations.
\newblock Mathematics of Computation. \textbf{43}, 1--19 (1984)

\bibitem{bfm}
Cristiani, E.: A fast marching method for {H}amilton-{J}acobi equations
  modeling monotone front propagation.
\newblock J Sci Comput \textbf{39}, 189--205 (2009)

\bibitem{dijkstra}
Dijkstra, E.W.: A note on two problems in connexion with graphs.
\newblock Numerische Mathematik \textbf{1}, 269--271 (1959)

\bibitem{mesh2d}
Engwirda, D. MESH2D - Automatic Mesh Generation (2011)

\bibitem{evans}
Evans, L.: Partial Differential Equations, \emph{Graduate Studies in
  Mathematics}, vol.~19, 2nd edn.
\newblock American Mathematical Society (2010)

\bibitem{falconeferretti}
Falcone, M., Ferretti, R.: Discrete time high-order schemes for viscosity
  solutions of {H}amilton-{J}acobi-{B}ellman equations.
\newblock Numer. Math \textbf{67}, 315--344 (1994)

\bibitem{oumboundarytracking}
Frew, E.: Combining area patrol, perimeter surveillance, and target tracking
  using ordered upwind methods.
\newblock In: IEEE International Conference on Robotics and Automation, Kobe,
  Japan, pp. 3123--3128 (2009)

\bibitem{gonrof}
Gonzalez, R., Rofman, E.: On deterministic control problems: An approximation
  procedure for the optimal cost {I}. the stationary problem.
\newblock SIAM J. Control and Optimization \textbf{23}, 242--266 (1985)

\bibitem{oumfoldsinstructuralgeology}
Hjelle, O., Petersen, A.: A {H}amilton-{J}acobi framework for modeling folds in
  structural geology.
\newblock Math Geosci \textbf{43}, 741--761 (2011)

\bibitem{fmmtriangle}
Kimmel, R., Sethian, J.: Computing geodesic paths on manifolds.
\newblock Proc. Natl. Acad. Sci. USA \textbf{95}, 8431--8435 (1998)

\bibitem{kumarvlad}
Kumar, A., Vladimirsky, A.: An efficient method for multiobjective optimal
  control and optimal control subject to integral constraints.
\newblock Journal of Computational Mathematics \textbf{28}, 517--551 (2010)

\bibitem{monneau}
Monneau, R.: Introduction to the {F}ast {M}arching {M}ethod.
\newblock Tech. rep., Centre International de Math{\'e}matiques Pures et
  Appliqu{\'e}s (2010)

\bibitem{oum}
Sethian, J., Vladimirsky, A.: Ordered upwind methods for static
  {H}amilton-{J}acobi equations: Theory and algorithms.
\newblock SIAM J. Numer. Anal. \textbf{41}, 325--363 (2003)

\bibitem{shumpaper}
Shum, A., Morris K.A., Khajepour, A., Direction-dependent optimal path planning for autonomous vehicles. \newblock Robotics and Autonomous Systems, \textbf{70}, 202-214 (2015)

\bibitem{souganidis}
Souganidis, P.: Approximation scheme for viscosity solutions of {H}amilton-{J}acobi equations.
\newblock Journal of Differential Equations. \textbf{59}, 1--43 (1985)

\bibitem{thepathtoefficiency}
Valero-G\'{o}mez, A., G\'{o}mez, J.V., Garrido, S., Moreno, L.: The path to efficiency: Fast Marching Method for safer, more efficient mobile robot trajectories.
\newblock IEEE Robotics and Automation Magazine. \textbf{20}(4), 111--120 (2013)

\bibitem{vladthesis}
Vladimirsky, A.: Fast methods for static {H}amilton-{J}acobi partial
  differential equations.
\newblock Ph.D. thesis, University of Califoria, Berkeley (2001)

\end{thebibliography}

\end{document}